\documentclass[reqno,11pt]{article}

\usepackage{amssymb,amsmath}

\usepackage{epsfig}
\usepackage{graphicx,color}

\usepackage[latin1]{inputenc}

\newtheorem{theo}{Theorem}
\newtheorem{cor}[theo]{Corollary} \newtheorem{lem}[theo]{Lemma}
\newtheorem{prop}[theo]{Proposition}
\newtheorem{defn}[theo]{Definition}
\newtheorem{rem}{Remark}
\newtheorem{ex}{Example}

\def\E{\mathbb{E}}

\def\C{\mathbb{C}}

\def\lip{{\rm Lip}}

 \def\d{\ \mathrm{d}}

\def\phi{\varphi}

\newcommand {\nn}{\nonumber}

\newcommand {\noi}{\noindent}

\def\mathsf{\bf}
\def\N{\mathbb{N}}

\def\R{\mathbb{R}}
\def\Z{\mathbb{Z}}

\def\d{\mathrm d}

\def\E{\mathrm E}
\def\P{\mathrm P}

\def\text{\mbox}

\def\1{{\bf 1}}

\newcommand\beqn{\begin{displaymath}}  
\newcommand\eeqn{\end{displaymath}}

\begin{document}

\title{
A nonlinear model for long memory
conditional heteroscedasticity}
\author{Paul  Doukhan$^1$, \ Ieva Grublyt\.e$^{1,2}$ \ and \ Donatas Surgailis$^2$  }
\date{\today \\  \small
\vskip.2cm
$^1$ Universit\'e de Cergy Pontoise, D\'epartament de Math\'ematiques, 95302 Cedex, France \\
$^2$ Vilnius University, Institute of Mathematics and Informatics, 08663 Vilnius, Lithuania}
\maketitle

\begin{quote}

{\bf Abstract.} We discuss a class of conditionally heteroscedastic time series models satisfying the equation
$r_t= \zeta_t \sigma_t$, where
$\zeta_t$ are standardized i.i.d. r.v.'s and  the conditional standard deviation
$\sigma_t$
is a nonlinear function $Q$ of inhomogeneous linear combination
of past values $r_s, s<t$ with coefficients $b_j$.
The existence of stationary solution $r_t$
with finite $p$th moment, $0< p < \infty $
is obtained under some conditions on $Q,  b_j$ and the $p$th moment
of $\zeta_0$.  Weak dependence
properties of $r_t$ are studied,
including the invariance principle for partial sums of Lipschitz functions of $r_t$.
In the case when $Q$
is the square root of a quadratic polynomial,
we prove that $r_t$ can exhibit a leverage effect and long memory,
in the sense that
the squared process $r^2_t$ has long memory autocorrelation and its  
normalized partial sums process
converges to a fractional Brownian motion.

\medskip


\end{quote}

\section{Introduction}

A stationary time series $\{r_t, t \in \Z \}$ is said {\it conditionally heteroscedastic}
if its conditional variance $\sigma^2_t = {\rm Var} [r_t |r_s, s< t] $ is a non-constant random process.
A class of  conditionally heteroscedastic
ARCH-type processes is defined from a standardized i.i.d. sequence $\{ \zeta_t, t \in \Z\}$ as solutions
of stochastic equation
\begin{equation}\label{genform}
r_t = \zeta_t \sigma_t, \qquad \sigma_t = V(r_s, s < t),
\end{equation}
where $V(x_1,x_2, \cdots)$ is some function of  $x_1,x_2, \cdots$. \\
The ARCH$(\infty)$ model 
corresponds to $V(x_1,x_2, \cdots) = \big(a + \sum_{j=1}^\infty b_j  x^2_j\big)^{1/2}$, or
\begin{equation}\label{ARCHinf}
\sigma^2_t = a + \sum_{j=1}^\infty b_j r^2_{t-j},
\end{equation}
where $a\ge 0, b_j\ge 0$ are coefficients. \\
The ARCH$(\infty)$ model includes
the well-known ARCH$(p)$ and GARCH$(p,q)$ models of Engle \cite{eng1982} and Bollerslev \cite{bol1986}.
However, despite their tremendous success, the GARCH models are not able to capture some
empirical features of asset returns, in particularly, the asymmetric or leverage effect discovered
by Black \cite{bla1976}, and the long memory decay in autocorrelation of squares $\{r^2_t\}$.
Giraitis and Surgailis \cite{gir2002} proved that the squared stationary solution
of the ARCH($\infty$) model in \eqref{ARCHinf} with $a>0$  always has short memory,
in the sense that $\sum_{j=0}^\infty {\rm Cov}(r^2_0, r^2_j) < \infty $. (However, for  integrated
ARCH($\infty$) models with $\sum_{j=1}^\infty b_j =1,  b_j \ge 0 $ and $a=0$ the situation is different;
see \cite{gir2014}.)

The above shortcomings of  the ARCH($\infty$) model motivated numerous studies proposing  alternative forms of the conditional variance
and  the function $V(x_1,x_2, \cdots)$  in \eqref{genform}. In particular, stochastic volatility models
can display both long memory and leverage except that in their case, the conditional variance is not a function of
$r_s, s< t$ alone and therefore it is more difficult to estimate from real data in comparison with the ARCH models;
see Shephard and Andersen \cite{she2009}.
Sentana \cite{sen1995} discussed a class of Quadratic ARCH (QARCH) models
with $\sigma^2_t $ 
being a general quadratic form in lagged variables $r_{t-1}, \cdots, $  $ r_{t-p}$. Sentana's specification of $\sigma^2_t$
encompasses a variety of ARCH models including
the asymmetric ARCH model of Engle \cite{eng1990} and the `linear standard deviation'
model of Robinson \cite{rob1991}. The limiting case (when $p=\infty$)
of the last model is the LARCH model discussed in \cite{gir2000} (see also \cite{gir2002}, \cite{ber2003}, \cite{gir2004}, \cite{tru2014}) and
corresponding to $V(x_1,x_2, \cdots) = a + \sum_{j=1}^\infty b_j  x_j,$ or
\begin{equation}\label{LARCHform}
\sigma_t = a + \sum_{j=1}^\infty b_j  r_{t-j},
\end{equation}
where  $a \in \R, b_j\in \R$ are real-valued coefficients. \cite{gir2000} proved
that the squared stationary solution  $\{r^2_t\}$ of the LARCH model with $b_j$ decaying as $j^{d-1}, 0< d< 1/2 $
may have long memory autocorrelations. The leverage effect in the LARCH
model was discussed in detail in \cite{gir2004}. On the other hand, volatility $\sigma_t $ (\ref{LARCHform})
of the LARCH model may assume negative values, lacking some
of the usual volatility interpretation.

The present paper discusses a class of conditionally heteroscedastic  models  \eqref{genform} with $V$ of the form
\begin{equation}\label{Qform}
V(x_1,x_2, \cdots)  =  Q(a + \sum_{j=1}^\infty b_j  x_j),
\end{equation}
where $Q(x), x \in \R$ is a (nonlinear) function of a single real variable $x\in \R$ which may be separated from zero by a positive
constant: $Q(x) \ge c >0, \, x \in \R$.  Linear $Q(x) = x$ corresponds to the LARCH model (\ref{LARCHform}).
Probably, the most interesting nonlinear case of $Q$ in (\ref{Qform}) is
\begin{equation*}\label{Qform1}
Q(x) = \sqrt{c^2 +  x^2},
\end{equation*}
where $c \ge 0$ is a parameter. In the latter case, the model is described by equations
\begin{eqnarray}
r_t&=&\zeta_t \sigma_t, \qquad
\sigma_t\ =\ \sqrt{c^2 + \Big(a + \sum_{s< t} b_{t-s} r_s\Big)^2}.
\label{newQ}
\end{eqnarray}
Note $\sigma_t \ge c \ge 0$ in (\ref{newQ}) is nonnegative and separated from $0$ if $c>0$.
Particular cases of volatility forms
in \eqref{newQ}  are:
\begin{eqnarray}
\sigma_t&=&\sqrt{c^2 + (a + br_{t-1})^2}   \qquad \text{(Engle's \cite{eng1990} asymmetric ARCH(1))}, \label{M1} \\
\sigma_t&=&\sqrt{c^2 + \big(a + \frac{b}{p}\sum_{j=1}^p r_{t-j}\big)^2}, \label{M2} \\
\sigma_t&=&\big|a + \sum_{j=1}^\infty b_j r_{t-j}\big|  \qquad \text{($Q(x)= |x|$)}, \label{M4}\\
\sigma_t&=&\sqrt{c^2    + (a + b((1-L)^{-d}-1)r_t)^2}. \label{M5}
\end{eqnarray}
\noi In  \eqref{M1}-\eqref{M5},  $a,b,c$ are real parameters, $p \ge 1$ an integer, $Lx_t =  x_{t-1}$ is the backward shift, and  $(1-L)^{-d} x_t = \sum_{j=0}^\infty \phi_j x_{t-j},  \, \phi_j =  \Gamma (d +j)/\Gamma(d) \Gamma(j+1),
\phi_0 =1 $ is the fractional
integration operator,  $0< d < 1/2$. The squared volatility (conditional variance) $\sigma^2_t$
in \eqref{M1}-\eqref{M5} and \eqref{newQ} is a quadratic form in lagged returns $r_{t-1}, r_{t-2}, \cdots $ and hence
represent particular cases of 
Sentana's \cite{sen1995} Quadratic ARCH (QARCH) model with $p = \infty $.
It should be noted, however, that the first two conditional moments
do not determine the unconditional distribution.  Particularly,
\eqref{genform} with  \eqref{newQ} generally  is a different process from  Sentana's \cite{sen1995} QARCH process, the latter being
defined as a solution of a linear random-coefficient equation for $\{r_t\}$ in contrast to the nonlinear equation in \eqref{genform}.
See also Example 2 below.

Let us describe the main results of this paper. Section 2 obtains sufficient
conditions on $Q, b_j $ and
$|\mu|_p := \E |\zeta_0|^p$ for the existence of  stationary solution of \eqref{genform}-\eqref{Qform}
with finite moment $\E |r_t|^p < \infty, \, p >0 $. We use the fact that the above equations can be
reduced to the `nonlinear moving-average' equation
\begin{equation*}
X_t = \sum_{s <t} b_{t-s} \zeta_s Q(a + X_s)
\end{equation*}
for linear form $X_t := \sum_{s<t} b_{t-s} r_s $ in \eqref{Qform}, and vice-versa.
Section 3 aims at providing weak dependence properties of  \eqref{genform} with $V$ in \eqref{Qform}, in particular,
the invariance principle for Lipschitz functions of $\{r_t\}$ and $\{X_t\}$,  under
the assumption that  $b_j$ are summable and decay as $j^{-\gamma}$ with $\gamma >1$.
Section 4 discusses long memory property of the quadratic model in  \eqref{newQ}.
For $b_j \sim \beta j^{d-1}, \, j \to \infty, \, 0< d < 1/2$ as in \eqref{M5},
we prove
that the squared process $\{r^2_t\} $
has long memory autocorrelations and its normalized partial sums process tend
to a fractional Brownian motion with Hurst parameter $H = d + 1/2 $ (Theorem \ref{long}).
Finally Section 5 establishes the leverage effect in spirit of \cite{gir2004}, viz.,
the fact that the `leverage function' $h_j := {\rm Cov}(\sigma^2_t, r_{t-j}), j \ge 1 $ of model \eqref{newQ}
takes negative values 
provided the coefficients $a$ and $b_j$ have opposite signs.
All proofs are collected in Section 6 (Appendix).

\smallskip

{\it Notation.} In what follows, $C, C(\dots)$ denote generic constants, possibly dependent
on the variables in brackets,
which may be different at different locations. $a_t \sim b_t \, (t \to \infty)$ is equivalent to
$\lim_{t\to \infty} a_t/b_t =1 $.

\section{Stationary solution}

This section discusses the existence of a stationary solution of  \eqref{genform} with
$V$ of \eqref{Qform}, viz.,
\begin{equation}\label{genformQ}
r_t = \zeta_t Q\big(a + \sum_{s<t} b_{t-s} r_s\big), \qquad t \in \Z.
\end{equation}
Denote
\begin{equation}\label{Xdef}
X_t \ := \ \sum_{s< t} b_{t-s} r_s.
\end{equation}
Then $r_t$ in \eqref{genformQ} can be written as  $r_t = \zeta_t Q(a + X_t)$ where \eqref{Xdef} formally satisfies
the following equation:
\begin{equation}\label{Xform}
X_t\ = \  \sum_{s<t}  b_{t-s} \zeta_s Q(a + X_s).
\end{equation}
Below we give
rigorous definitions of solutions of  \eqref{genformQ} and  \eqref{Xform} and a statement  (Proposition \ref{Xreq}) justifying
\eqref{Xform} and the equivalence of \eqref{genformQ} and  \eqref{Xform}.

In this section we consider a general case of \eqref{genformQ}-\eqref{Xform}
when the innovations may have infinite variance. More precisely, we assume that
$\{\zeta_t, t\in \Z\}$ are i.i.d. r.v.'s with finite
moment $|\mu|_p :=  \E |\zeta_t|^p < \infty, \, p>0$.
We use the following moment inequality.

\begin{prop} \label{Yp} Let
$\{Y_j, j \ge 1\}$ be a sequence of r.v.'s such that $\E |Y_j|^p < \infty $ for some
$p >0$ and the sum on the r.h.s. of \eqref{rosen} converges.
If $p>1 $ we additionally assume that $\{Y_j\}$ is a martingale difference sequence:
$\E [Y_j |Y_1, \cdots, Y_{j-1}] = 0, \, j=2,3, \cdots $. Then
there exists a constant $K_p$ depending only on $p$ and such that
\begin{equation}\label{rosen}
\E \big|\sum_{j=1}^\infty Y_j\big|^p \ \le \ K_p \begin{cases}
\sum_{j=1}^\infty \E |Y_j|^p, &0< p \le 2, \\
\big(\sum_{j=1}^\infty (\E |Y_j|^p)^{2/p}\big)^{p/2}, &p > 2.
\end{cases}
\end{equation}

\end{prop}

\begin{rem} \label{rosenC} {\rm In the sequel, we shall refer to $K_p$ in \eqref{rosen} as the Rosenthal constant.
For $0< p \le 1$ and $p=2$, inequality \eqref{rosen} holds with $K_p = 1$, and
for $1 < p < 2 $,
it is known as von Bahr and Ess\'een inequality, see \cite{von1965}, which
holds with $K_p = 2$. For $p> 2$,  inequality \eqref{rosen}
is a consequence of the Burkholder and Rosenthal inequality (see \cite{bur1973}, \cite{ros1970}, also \cite{gir2012}, Lemma 2.5.2).
Os\c ekowski \cite{ose2012} proved that $K^{1/p}_p
\le 2^{(3/2) + (1/p)} (\frac{p}{4} + 1)^{1/p} \big(1 + \frac{p}{\log (p/2)}\big)$, in particular,
$K_4^{1/4} \le  27.083
$.
See also \cite{hit1990}. }
\end{rem}

\smallskip

Let us give some formal  definitions. Let ${\cal F}_t = \sigma(\zeta_s, s\le t), t \in \Z$ be the sigma-field generated
by $\zeta_s, s\le t$. A random process $\{ u_t, t \in \Z\}$ is called {\it adapted} (respectively,  {\it predictable})
if $u_t$ is ${\cal F}_t$-measurable for each $t \in \Z $ (respectively,  $u_t$ is ${\cal F}_{t-1}$-measurable for each $t \in \Z $). Define
\begin{equation}\label{Bp}
B_p := \begin{cases}
\sum_{j=1}^\infty |b_j|^p, &0<p < 2, \\
\big(\sum_{j=1}^\infty b_j^2\big)^{p/2}, &p\ge 2.
\end{cases}
\end{equation}

\begin{defn} Let $p>0$ be arbitrary.

\smallskip

\noi (i) By $L^p$-solution of   \eqref{genformQ} we mean  an adapted process $\{r_t, t \in \Z\}$
with $\E |r_t|^p < \infty $ such that  for any $t \in \Z$ the series
$\sum_{s<t}  b_{t-s} r_s$ converges in $L^p$ and   \eqref{genformQ} holds.

\smallskip

\noi (ii) By $L^p$-solution of \eqref{Xform} we mean  an predictable process $\{X_t, t \in \Z\}$ with
$\E |X_t|^p < \infty $ such that  for any $t \in \Z$ the series
$\sum_{s<t}  b_{t-s} \zeta_s  Q(a+ X_s)$ converges in $L^p$ and \eqref{Xform}  holds.

\end{defn}

Let $Q(x), x \in \R $ 
be a Lipschitz function, i.e., there exists $\lip_Q >0$ such that
\begin{equation}\label{QLip}
|Q(x)-Q(y)| \le \lip_Q |x-y|, \qquad  x,y \in \R.
\end{equation}
Note \eqref{QLip} implies the bound
\begin{equation}\label{Qnel}
Q^2(x) \le c_1^2 + c_2^2 x^2, \qquad x \in \R,
\end{equation}
where $c_1\ge 0, \, c_2 \ge \lip_Q$ and $c_2 $ can be chosen arbitrarily close to $\lip_Q$;
in particular, \eqref{Qnel} holds with $c^2_2 = (1+ \epsilon^2)\lip^2_Q, c_1^2  = Q^2(0)(1 + \epsilon^{-2}), $ where $\epsilon>0$
is arbitrarily small.

\begin{prop} \label{Xreq} Let $Q$ be a measurable function satisfying  \eqref{Qnel} with some $c_1, c_2 \ge 0$ and $\{\zeta_t\}$ be an i.i.d. sequence with $|\mu|_p = \E|\zeta_0|^p < \infty$ and satisfying
$\E \zeta_0 = 0$ for $p> 1$. In addition, assume $B_p < \infty$.

\smallskip

\noi (i) Let $\{X_t \}$ be a stationary $L^p$-solution of \eqref{Xform}. Then $\{r_t  := \zeta_t Q(a + X_t)\}$ is a stationary
$L^p$-solution of \eqref{genformQ} and
\begin{equation}\label{rX}
\E |r_t|^p \ \le \ C(1+ \E |X_t|^p).
\end{equation}
Moreover, for $p>1 $, $\{r_t, {\cal F}_t, t \in \Z\}$ is a
martingale difference sequence with
\begin{equation}\label{genformvol}
\E [r_t|{\cal F}_{t-1}] = 0, \qquad \E [|r_t|^p| {\cal F}_{t-1}]  = |\mu|_p \big|Q(a + \sum_{s<t}  b_{t-s} r_s)|^p.
\end{equation}

\noi (ii) Let $\{r_t \}$ be a stationary $L^p$-solution of \eqref{genformQ}.
Then  $\{X_t\}$ in \eqref{Xdef} is a   stationary $L^p$-solution of \eqref{Xform} such that
\begin{equation}\label{Xr}
\E |X_t|^p \ \le \ C\E |r_t|^p.
\end{equation}
Moreover, for $p \ge 2$
\begin{equation}\label{covX}
\E [X_t X_0] \ = \ \E r^2_0 \sum_{s=1}^\infty b_{t+s} b_s, \qquad t = 0,1,\dots.
\end{equation}
\end{prop}

\begin{rem} \label{Lpsol} {\rm  Let $p\ge 2$ and $|\mu|_p < \infty$, then  by inequality \eqref{rosen},
$\{r_t\}$ being a stationary $L^p$-solution of \eqref{genformQ} is equivalent
to  $\{r_t\}$ being a stationary $L^2$-solution of \eqref{genformQ} with $\E |r_0|^p < \infty $.
Similarly, if $Q$ and $\{\zeta_t\}$  satisfy the conditions of Proposition \ref{Xreq} and $p\ge 2$,
then  $\{X_t\}$ being a stationary $L^p$-solution of \eqref{Xform} is equivalent
to  $\{X_t\}$ being a stationary $L^2$-solution of \eqref{Xform} with $\E |X_0|^p < \infty $.
}
\end{rem}

The following theorem obtains a sufficient condition in \eqref{cQB} for the existence of a stationary
$L^p$-solution of equations \eqref{genformQ} and \eqref{Xform}. Condition  \eqref{cQB} involves
the $p$th moment of innovations, the Lipschitz constant $\lip_Q$, the sum $B_p$ in \eqref{Bp} and the Rosenthal
constant $K_p$ in \eqref{rosen}. Part (ii) of Theorem \ref{Xexists} shows that for $p=2$, condition
\eqref{cQB} is close to optimal, being necessary in the case of quadratic $Q^2(x) = c_1^2 + c_2^2 x^2 $.

\begin{theo} \label{Xexists} Let the conditions of Proposition \ref{Xreq} be satisfied, $p>0$ is arbitrary.
In addition,
assume that $Q$ satisfies
the Lipschitz condition in \eqref{QLip}.

\medskip

\noi (i) Let
\begin{equation} \label{cQB}
K_p |\mu|_p \lip^p_Q B_p   < 1.
\end{equation}
Then there exists a unique stationary $L^p$-solution $\{X_t\}$ of \eqref{Xform} and
\begin{equation}\label{X2mom}
\E |X_t|^p \ \le \frac{C(p,Q) |\mu|_p  B_p }{1 - K_p |\mu|_p \lip^p_Q B_p},
\end{equation}
where $C(p,Q) < \infty $ depends only on $p$ and $c_1, c_2 $ in \eqref{Qnel}.

\smallskip

\noi (ii) Assume, in addition, that
$Q^2(x) = c_1^2 + c_2^2 x^2$, where $c_i \ge 0, i=1,2$, and $\mu_2 = \E \zeta^2_0 =1 $.
Then $c^2_2 B_2 < 1$ is a necessary and sufficient condition for the existence
of  a stationary $L^2$-solution $\{X_t\}$ of \eqref{Xform} with $a \neq 0$.
\end{theo}

\begin{rem}{\rm
Condition \eqref{cQB} agrees with the contraction condition for the operator
defined by the r.h.s. of \eqref{Xform} and acting in a suitable space of predictable processes
with values in $L^p$.
For the LARCH model,
explicit conditions for finiteness of the $p$th moment were obtained in
\cite{gir2000}, \cite{gir2004} using a specific diagram approach for multiple Volterra series.
For larger values of $p>2 $, condition \eqref{cQB} is preferable to the corresponding condition
\begin{equation}\label{larchp}
(2^p - p -1)^{1/2} |\mu|_p^{1/p} B_p^{1/p} < 1, \qquad p =2,4,6, \cdots,
\end{equation}
in (\cite{gir2000}, (2.12)) for the LARCH model, since the coefficient $ (2^p - p -1)^{1/2}$ grows exponentially
with $p$ in contrast to the bound on $K_p^{1/p} $ in Remark \ref{rosenC}.
See also (\cite{gru2014}, sec.~4.3).  On the other hand for $p=4$ \eqref{larchp}
becomes $\sqrt{11} |\mu|_4^{1/4} B_2^{1/2} < 1 $ while \eqref{cQB}
is satisfied if $K_4^{1/4} |\mu|^{1/4}_4 B_2^{1/2} \le  27.083 |\mu|^{1/4}_4 B_2^{1/2} < 1 $,
see   Remark \ref{rosenC}, which is worse than   \eqref{larchp}.

}
\end{rem}

\begin{ex} [The LARCH model]\label{exlarch} {\rm Let $Q(x) = x$ and $\{\zeta_t\}$ be a standardized
i.i.d. sequence with zero mean and unit variance.
Then 
\eqref{Xform} becomes the bilinear equation
\begin{equation}\label{Lform}
X_t\ = \  \sum_{s<t}  b_{t-s} \zeta_s (a + X_s).
\end{equation}
The corresponding conditionally heteroscedastic process $\{r_t = \zeta_t (a + X_t)\}$ in Proposition \ref{Xreq}(i) is  the LARCH model discussed
in \cite{gir2000}, \cite{gir2004} and elsewhere. As shown in (\cite{gir2000}, Thm.2.1),
equation
\eqref{Lform} admits a covariance stationary predictable solution if and only if $B_2 = \sum_{j=1}^\infty b^2_j < 1$.
Note the last result agrees  with Theorem \ref{Xexists} (ii).
A crucial role in the study of the LARCH model is played by the fact that its solution
can be written in terms of the convergent orthogonal  Volterra series
$$
X_t = a \sum_{k=1}^\infty \sum_{s_k < \cdots < s_1 < t} b_{t-s_1} \cdots b_{s_{k-1}-s_k} \zeta_{s_1} \cdots \zeta_{s_k}.
$$
Except for $Q(x) = x$, in other cases of \eqref{Xform} including the QARCH model in \eqref{newQ},
Volterra series expansions are unknown and their usefulness is doubtful.

 }
\end{ex}

\begin{ex} [Asymmetric ARCH(1)]  \label{Qarch}  {\rm Consider the model \eqref{genform} with
$\sigma_t $ in \eqref{M1}, viz.
\begin{equation}
r_t = \zeta_t \big(c^2 + (a + br_{t-1})^2\big)^{1/2},    \label{M11}
\end{equation}
where $\{\zeta_t\} $ are standardized i.i.d. r.v.'s. By Theorem \ref{Xexists} (ii), equation \eqref{M11} has a unique stationary
solution with finite variance $\E r^2_t = (a^2+c^2)/(1-b^2) $  if and only if $b^2 < 1$. \\
In parallel, consider the random-coefficient AR(1) equation
\begin{equation}
\widetilde r_t = \kappa \varepsilon_t  + b \eta_t \widetilde r_{t-1},    \label{M12}
\end{equation}
where
$\{ (\varepsilon_t, \eta_t)\} $ are i.i.d. random vectors with zero mean
$\E \varepsilon_t = \E \eta_t = 0 $ and unit variances $\E[ \varepsilon^2_t] = \E [\eta^2_t] =1$
and $\kappa, b$ are real coefficients.
As shown in Sentana \cite{sen1995} (see also Surgailis \cite{sur2008}), equation \eqref{M12} has a stationary solution
with finite variance under the same condition $b^2 < 1 $  as \eqref{M11}. Moreover, if the coefficients
$\kappa $ and $\rho := \E [\varepsilon_t \eta_t] \in [-1,1] $ in \eqref{M12}
are related to the coefficients $a, c $ in \eqref{M11} as
\begin{equation} \label{kappa}
\kappa \rho = a, \qquad \kappa^2 = a^2 + c^2,
\end{equation}
then the processes in \eqref{M11} and  \eqref{M12} have the same volatility forms since
\begin{eqnarray*}
\widetilde \sigma^2_t\ :=\ \E [\widetilde r^2_t| \widetilde r_s, s<t]
&=&\kappa^2 + 2 \kappa b \rho + b^2 \widetilde r^2_{t-1} \\
&=&c^2 + (a + b \widetilde r_{t-1})^2
\end{eqnarray*}
agrees with the corresponding expression  $\sigma^2_t = c^2 + (a + b r_{t-1})^2$  in the case  of \eqref{M1}.

A natural question is whether the above stationary solutions  $\{r_t\}$ and $\{\widetilde r_t\}$ of
\eqref{M11} and \eqref{M12}, with parameters related as in \eqref{kappa},
have the same (unconditional) finite-dimensional distributions? As shown in (\cite{sur2008}, Corollary 2.1),
the answer is positive in the case when $\{\zeta_t\}$ and $\{(\varepsilon_t, \eta_t)\} $ are
Gaussian sequences. However, the conditionally Gaussian case seems to be the only exception and
in general the processes
$\{r_t\}$ and $\{\widetilde r_t\}$
have different distributions.
This can be seen
by considering the 3rd conditional moment of \eqref{M11}
\begin{equation} \label{M1c}
\E [r^3_t|r_{t-1}] = \mu_3 \big(c^2 + (a + br_{t-1})^2\big)^{3/2}
\end{equation}
which is an irrational function of $r_{t-1}$ (unless $\mu_3= \E \zeta^3_0 = 0$ or $b=0$),
while a similar moment of \eqref{M12}
\begin{equation} \label{M2c}
\E [\widetilde r^3_t|\widetilde r_{t-1}] = \kappa^3 \nu_{3,0}  + 3 b \kappa^2 \nu_{2,1} \widetilde r_{t-1}
+ 3 b^2 \kappa \nu_{1,2} \widetilde r^2_{t-1} + b^3 \nu_{0,3} \widetilde r^3_{t-1}
\end{equation}
is a cubic polynomial in $\widetilde r_{t-1}$, where $\nu_{i,j} := \E [\varepsilon^i_0 \eta^{j}_0].$ Moreover,
\eqref{M1c} has a constant sign independent of $r_{t-1}$ while
the sign of the cubic polynomial in
\eqref{M2c} changes with $\widetilde r_t $ ranging from $\infty $ to $-\infty $
if the leading coefficient $b^3 \nu_{0,3} \ne 0$.

Using the last observation we can prove that the bivariate distributions of $(r_t, r_{t-1})$ and $(\widetilde r_t, \widetilde r_{t-1})$
are different under general  conditions on the innovations and the parameters of the two equations. The argument is as follows.
Let  $b > 0, c >0, \mu_3 >0, \nu_{0,3} = \E \eta_0^3 >0$.
Assume that
$\zeta_0 $ has a bounded strictly positive density function $0<f(x) < C, x \in \R $ and
$(\varepsilon_0,\eta_0)$ has a bounded strictly positive density function $0< g(x,y) <C, (x,y) \in \R^2$.
The above assumptions imply that the distributions of $r_t $ and $\widetilde r_t$ have infinite support. Indeed,
by \eqref{M11} and the above assumptions we have that  $\P(r_t >K ) = \int_{\R} \P( c^2 + (a + b r_{t-1})^2 > (K/y)^2) f(y) \d y  >0 $ for any $K>0$ since $\lim_{y \to \infty} \P( c^2 + (a + b r_{t-1})^2 > (K/y)^2) = 1 $.  Similarly,
 $\P(\widetilde r_t >K ) = \int_{\R^2} \P( \widetilde r_{t-1} > (K- \kappa x)/b y ) g(x,y) \d x \d y  >0 $ and
$ \P(r_t < - K )> 0,   \P(\widetilde r_t < - K )  >0$ for any $K >0$. Since
$h(x) := \mu_3 \big(c^2 + (a + bx)^2\big)^{3/2} \ge 1 $ for all $ |x|> K$ and any sufficiently large $K >0$, from
\eqref{M1c} we obtain that for any $K>0$
\begin{eqnarray}
\E r^3_t \1( r_{t-1} >K)&=&\E h(r_{t-1}) \1(r_{t-1} >K) > 0   \quad  \text{and} \nn \\
\E r^3_t \1( r_{t-1} <- K)&=&\E h(r_{t-1}) \1(r_{t-1}<-K) > 0. \label{h1}
\end{eqnarray}
On the other hand, since   $\widetilde h(x) := \kappa^3 \nu_{3,0}  + 3 b \kappa^2 \nu_{2,1} x
+ 3 b^2 \kappa \nu_{1,2} x^2 + b^3 \nu_{0,3} x^3  \ge 1 $ for $x > K $ and $\widetilde h(x) \le - 1 $ for $x < - K $
and $K$ large enough, from
\eqref{M2c} we obtain that for all sufficiently large $K>0$
\begin{eqnarray}
\E \widetilde r^3_t \1(\widetilde r_{t-1} >K)&=&\E \widetilde h(\widetilde r_{t-1}) \1(\widetilde r_{t-1} >K) > 0   \quad  \text{and} \nn \\
\E \widetilde r^3_t \1(\widetilde r_{t-1} <- K)&=&\E \widetilde h(\widetilde r_{t-1}) \1(\widetilde r_{t-1}<-K) < 0. \label{h2}
\end{eqnarray}
Clearly, \eqref{h1} and \eqref{h2} imply that   the bivariate distributions of $(r_t, r_{t-1})$ and $(\widetilde r_t, \widetilde r_{t-1})$ are
different under the stated assumptions. 

\smallskip

For models \eqref{M11} and \eqref{M12}, we can explicitly compute covariances
$\rho(t)= {\rm cov}( r_t^2, r_0^2),  \widetilde \rho(t)= {\rm cov}(\widetilde r_t^2, \widetilde r_0^2)$ and some
other joint moment functions, as follows.

Let $\mu_3 = \E \zeta^3_0 = 0, \mu_4 = \E \zeta^4_0 < \infty $ and
$m_2 := \E r^2_0, \,  m_3(t) :=\E r_t^2 r_0,\, m_4(t) :=\E r_t^2 r_0^2, t \ge 0$.
Then
\begin{eqnarray}
m_2&=&(a^2+c^2)/(1-b^2), \qquad m_3(0) \ = \ 0,  \nn \\
m_3(1)&=&\E [((a^2+c^2) +2 ab r_{0} +b^2 r_{0}^2) r_0]
=2abm_2 + b^2 m_3(0) \ = \ 2abm_2, \nn  \\
m_3(t)&=&\E [((a^2+c^2) +2 ab r_{t-1} +b^2 r_{t-1}^2) r_0]
\ = \ b^2 m_3(t-1) \ = \cdots \ = \ b^{2(t-1)} m_3(1) \nn  \\
&=&\frac{2ab(a^2+c^2)}{1-b^2} b^{2(t-1)} , \quad t\ge 1. \label{m3}
\end{eqnarray}
Similarly,
\begin{eqnarray*}
m_4(0)&=&\mu_4\E [((a^2+c^2) +2 ab r_{-1} +b^2 r_{-1}^2)^2] \\
&=&\mu_4\{(a^2+c^2)^2 + (2ab)^2m_2 + b^4 m_4(0)
+ 2b^2 (a^2 +c^2)m_2\}, \\
m_4(t)&=&
\E [((a^2+c^2) +2 ab r_{t-1} +b^2 r_{t-1}^2) r^2_0]
=(a^2+c^2)m_2 + b^2 m_4(t-1), \quad t\ge 1
\end{eqnarray*}
resulting in
\begin{eqnarray}
m_4(0)&=&\frac{\mu_4((a^2+c^2)^2 +((2ab)^2+2(a^2+c^2) b^2)m_2)}{1- \mu_4b^4}, \label{m41}\\
m_4(t)&=&m_2(a^2+c^2)\cdot \frac{1-b^{2t}}{1- b^2} + b^{2t} m_4(0), \quad t \ge 1, \nn
\end{eqnarray}
and
\begin{eqnarray}\label{rhot}
\rho(t)&=&(m_4(0)- m_2^2)b^{2t}, \qquad t\ge 0.
\end{eqnarray}
In a similar way, when the distribution of $\zeta_0$ is symmetric one can write recursive linear equations
for joint even moments $\E[r^{2p}(0) r^{2p}(t)]$ of arbitrary order  $p = 1,2,\dots $
involving $ \E[r^{2l}(0) r^{2p}(t)], 1\le l \le p-1$ and $m_{2k}(0) = \E [r^{2k}(0)], 1\le k \le 2p$. These equations
can be explicitly solved in terms of $a,b,c$ and $\mu_{2k}, 1\le k \le 2p$.

A similar approach can be applied to find joint moments of
the random-coefficient AR(1) process in \eqref{M12}, with the difference that
symmetry of $(\varepsilon_0,\eta_0)$ is not needed.
Let  $\widetilde m_2:= \E \widetilde{r}_t^2$, $\widetilde{m}_3(t) :=\E [\widetilde{r}_t^2 \widetilde{r}_0]$, $\widetilde{m}_4(t) :=\E [\widetilde{r}_t^2 \widetilde{r}_0^2]$ and $\widetilde{\rho}(t) := {\rm Cov}( \widetilde{r}_t^2, \widetilde{r}_0^2), \, \nu_{i,i} := \E[\varepsilon_0^i \eta_0^j]$. Then
\begin{eqnarray*}
\widetilde m_2&=&\kappa^2/(1-b^2), \\
\widetilde{m}_3(0)&=&\E[(\kappa \varepsilon_0  + b \eta_0 \widetilde r_{-1})^3]
\ = \ \kappa^3 \nu_{3,0} + 3 \kappa b^2 \nu_{1,2} \widetilde m_2 + b^3 \nu_{0,3} \widetilde{m}_3(0), \\
\widetilde{m}_3(1)&=&
\E[ (\kappa +2 \kappa \rho b \widetilde{r}_{0} +b^2 \widetilde{r}_{0}^2) \widetilde{r}_0] \ = \
2 \kappa \rho b \widetilde m_2 +  b^2 \widetilde m_3(0), \\
\widetilde{m}_3(t)&=&
\E [(\kappa ^2 + 2 \kappa \rho b \tilde{r}_{t-1} + b^2 \tilde{r}_{t-1}^2) \tilde{r}_0] \ = \
b^2 \widetilde{m}_3(t-1) \ = \cdots  = \ b^{2(t-1)} \widetilde{m}_3(1), \quad t \ge 2,
\end{eqnarray*}
and
\begin{eqnarray*}
\widetilde m_4(0)&=&\E[(\kappa \varepsilon_0  + b \eta_0 \widetilde r_{-1})^4] \\
&=&\kappa^4 \nu_{4,0} +
6 \kappa^2 b^2 \nu_{2,2} \widetilde m_2+ 4\kappa b^3 \nu_{1,3} \widetilde m_3(0) + b^4 \nu_{0,4} \widetilde m_4(0),  \\
\widetilde m_4(1)&=&\E[(\kappa \varepsilon_t  + b \eta_t \widetilde r_{0})^2 \widetilde r^2_0]  =  \kappa^2 \widetilde m_2 + 2 \kappa  \rho b \widetilde m_3(0) + b^2 \widetilde {m}_4(0)\nn \\
\widetilde m_4(t)&=&\E[(\kappa \varepsilon_t  + b \eta_t \widetilde r_{t-1})^2 \widetilde r^2_0]
\ = \ \kappa^2 \widetilde m_2 + b^2 \widetilde {m}_4(t-1), \qquad t\ge 2,
\end{eqnarray*}
leading to
\begin{eqnarray}
\widetilde m_3(0)&=&\frac{\kappa^3 \nu_{3,0} + 3\kappa b^2 \nu_{1,2} \widetilde m_2 }{1- \nu_{0,3} b^3}, \nn \\
\widetilde m_3(t)&=& b^{2(t-1)}(  2 \kappa \rho b \widetilde m_2 + b^2 \widetilde m_3(0) )\quad t\ge 1.\nn \\
\widetilde m_4(0)&=&\frac{\kappa^4 \nu_{4,0} + 6 \kappa^2 b^2 \nu_{2,2} \widetilde m_2+ 4\kappa b^3 \nu_{1,3} \widetilde m_3(0)}{1-\nu_{0,4} b^4}, \label{m42}\\
\widetilde m_4(t)&=& \widetilde m_2\kappa^2 \Big(\frac{1-b^{2t}}{1- b^2}\Big)  + b^{2t}( \widetilde {m}_4(0)+2 \kappa  \rho \widetilde m_3(0)/b ), \quad t \ge 1, \nn \nn
\end{eqnarray}
and
\begin{eqnarray*}
\widetilde{\rho}_4(t)&=&b^{2(t-1)} \widetilde{\rho}_4(1), \quad t \ge 1,\\
\widetilde{\rho}_4(1)&=&2 \rho \kappa b \widetilde m_3(0) +b^2( \widetilde m_4(0)- \widetilde m_2^2).
\end{eqnarray*}
Then if $\nu_{3,0} = \nu_{1,2} = 0$ we have $\widetilde m_3(0)= 0$ and
$\widetilde{\rho}_4(t)= ( \widetilde m_4(0)- \widetilde m_2^2) b^{2t}$; moreover,
$\widetilde m_2 = m_2$ in view of \eqref{kappa}. Then
$\widetilde{\rho}_4(t) = {\rho}_4(t)$ is equivalent to
$\widetilde m_4(0) = m_4(0)$, which follows from
\begin{equation}\label{m43}
\mu_4 = \nu_{0,4} = \nu_{4,0} \quad \text{and} \quad 6 \nu_{2,2} = \mu_4 (4\nu^2_{1,1} + 2),
\end{equation}
see \eqref{rhot}, \eqref{m41}, \eqref{m42}. Note that  \eqref{m43} hold for centered Gaussian distribution
$(\varepsilon_0, \eta_0)$ with unit variances $\E \varepsilon^2_0 = \E \eta^2_0 =1 $.

}
\end{ex}

\section{Weak dependence}

Various measures of weak dependence for stationary processes $\{y_t\} = \{y_t, t \in \Z\}$
have been introduced in the  literature, see  e.g. \cite{Ded2007}. Usually, the dependence between the present $(t\ge 0)$ and
the past $(t\le -n)$ values of $\{y_t\}$ is measured by some dependence coefficients
decaying to $0$ as $n \to \infty$. The decay rate of these coefficients
plays a crucial role in establishing many asymptotic results.
The classical problem is proving Donsker's invariance principle:
\begin{equation} \label{clt}
\frac1{\sqrt n}\sum_{t  =1}^{[n\tau]} (y -  \E y_t) \to \sigma B(\tau), \qquad \mbox{ in  the Skorohod space } \ D[0,1],
\end{equation}
where $B = \{B(\tau), \tau \in [0,1]\}$ is a standard Brownian motion. The above result is useful in change-point analysis
(Cs\"org\H{o} and Horv\'ath \cite{cso1997}), financial mathematics and many other areas. Further applications
of weak dependence coefficients include empirical  processes \cite{ded2007} and the asymptotic behavior of various
statistics, including the maximum likelihood estimators. See Ibragimov and Linnik \cite{ibr1971} and
the application to GARCH estimation in \cite{lin2009}. \\
The present Section discusses two measures of weak dependence - the projective weak dependence
coefficients of Wu \cite{wub2005} and the $\tau$-dependence coefficients introduced in
Dedecker and Prieur \cite{ded2004}, \cite{ded2005} - for stationary solutions $\{r_t\}, \{X_t\}$ of
equations  (\ref{genformQ}),  \eqref{Xform}. We show that
the decay rate of the above weak dependence coefficients is determined by the decay rate
of the moving average coefficients $b_j$.

\subsection{Projective weak dependence coefficients}

Let us introduce some notation. For  r.v. $\xi$, write $\|\xi\|_p := \E^{1/p} |\xi|^p, \, p \ge 1 $. Let $\{y_t, t\in \Z\}$
be a stationary causal Bernoulli shift in i.i.d. sequence $\{\zeta_t\}$, in other words,
$$
y_t =  f(\zeta_s, s \le t), \quad t \in \Z,
$$
where $f:  \R^\N \to  \R $ is a measurable function.
We also assume $\E y_0 = 0, \|y_0\|^2_2 = \E y^2_0 < \infty.$
Introduce the projective weak dependence coefficients
\begin{equation}\label{omega2}
\omega_2(i; \{y_t\})\ := \ \|f_{i}(\xi_0) - f_{i}(\xi'_0)\|_2,
\qquad \delta_2(i; \{y_t\})\ := \ \|f(\xi_i) - f(\xi'_i)\|_2,
\end{equation}
where $\xi_i := (\cdots, \zeta_{-1}, \zeta_0, \zeta_1, \cdots, \zeta_i), \ \xi'_i := (\cdots, \zeta_{-1}, \zeta'_0, \zeta_1, \cdots, \zeta_i)$,
$f_i(\xi_0) := \E [ f(\xi_i)| \xi_0] = \E [y_i |{\cal F}_0] $ is the conditional expectation
and $\{\zeta'_0, \zeta_t, t \in \Z\}$ are i.i.d. r.v.s.
Note the i.i.d. sequences $\xi $ and $\xi'_i$ coincide except for a single entry. Then
$\omega_2(i; \{y_t\})\le \delta_2(i; \{y_t\}), i \ge 0$ and condition
\begin{equation}\label{wu}
\sum_{k=0}^\infty \omega_2(k; \{y_t\}) \ < \infty
\end{equation}
guarantees the weak invariance principle in \eqref{clt} with $\sigma^2 := \sum_{j \in \Z} {\rm Cov}(y_0, y_j), $
see Wu \cite{wub2005}. The last series absolutely converges in view of \eqref{wu} and the bound in
(\cite{wub2005}, Thm.~1), implying $ |{\rm Cov}(y_0, y_j)| \le \sum_{k=0}^\infty \omega_2(k; \{y_t\}) \omega_2(k+j; \{y_t\}), j \ge 0$.  \\
Below, we verify
Wu's  condition \eqref{wu}
for $\{X_t\}, \{r_t\}$ in (\ref{Xform}),  (\ref{genformQ}). We assume that the coefficients $b_j$ decay as
$j^{-\gamma}$ with some $\gamma > 1 $, viz.,
\begin{equation}\label{bcond}
\exists \gamma >0, \  c>0: \qquad
|b_j| <  c j^{-\gamma}, \quad \forall \ j \ge 1.
\end{equation}

\begin{prop} \label{propomega} Let $Q$ satisfy the Lipschitz condition in \eqref{QLip}, $p \ge 1, \,
K_p |\mu|_p \lip^p_Q B_p < 1 $ (see \eqref{cQB}),
and $\{ X_t\}, \{r_t\} $ be stationary $L^p$-solutions of
\eqref{Xform},  \eqref{genformQ},  respectively. In addition, assume that
$b_j$ satisfy \eqref{bcond} with $\gamma > \max \{1/2, 1/p \}$. \\
Then
\begin{equation}\label{deltas}
\delta_p(k; \{X_t\}) = O(k^{-\gamma}) \qquad \text{and} \qquad
\delta_p(k; \{r_t\}) = O(k^{-\gamma}).
\end{equation}
\end{prop}
The next corollary follows from the above-mentioned result of Wu \cite{wub2005}, relations
$\delta_2(k; \{y_t\} )$   $  \le  C\delta_2(k; \{r_t\} ),  \,
\delta_2(k; \{z_t\} )  \le  C\delta_2(k; \{X_t\} )$ and
\eqref{deltas}.

\begin{cor} Let $\{y_t := h(r_t) \},  \{z_t := h(X_t) \}$, where $\{ X_t\}, \{r_t\} $ are as in
Proposition \ref{propomega} and $h:  \R \to \R$ is a Lipschitz function.  Then
\begin{equation}\label{clt1}
n^{-1/2} \sum_{t=1}^{[n\tau]}(y_t- \E y_t) \ \to_{D[0,1]} c_y B(\tau) \quad \text{and} \quad  n^{-1/2} \sum_{t=1}^{[n\tau]}
(z_t - \E z_t) \ \to_{D[0,1]} c_z B(\tau),
\end{equation}
where $B$ is a standard Brownian motion and
\begin{equation*}\label{covser}
c^2_y := \sum_{t \in \Z} {\rm Cov}( y_0,  y_t) < \infty, \qquad c^2_z := \sum_{t \in \Z} {\rm Cov}(z_0,  z_t) < \infty.
\end{equation*}
\end{cor}

\subsection{$\tau$-weak dependence coefficients}

Let $\{y_t, t \in \Z\}$ be a stationary process with $\|y_0\|_p < \infty, p \in [1, \infty]$.
Following Dedecker and Prieur \cite{ded2004}, \cite{ded2005}, we define the $\tau$-weak dependence coefficients
\begin{eqnarray*}
\tau_p(\{y_{j_i}\}_{1 \le i \le k })
&:=&\Big\|\sup_{f \in \Lambda_1(\R^k)} \big| \E \big[ f(y_{j_1}, \cdots, y_{j_k}) \big|  y_t, t \le 0 \big]
-  \E [f(y_{j_1}, \cdots, y_{j_k})] \big| \Big\|_p
\end{eqnarray*}
measuring the dependence between $\{y_t\}_{t \le 0}$ and $\{y_{j_i}\}_{1\le i \le k}, 0< j_1 < \cdots <  j_k $, and
\begin{eqnarray*}
\tau_p(n; \{y_j\})\  :=\  \sup_{k \ge 1}  k^{-1} \sup_{n \le j_1 < \cdots < j_k}  \tau_p(\{y_{j_i}\}_{1 \le i \le k }).
\end{eqnarray*}
Here, $\Lambda_1(\R^k)$ denotes the class of all Lipschitz functions $f: \R^k \to \R$ with
\begin{equation*}
|f(x_1, \cdots, x_k)-  f(y_1, \cdots, y_k)| \le \sum_{i=1}^k |x_i-y_i| \quad \text{for any} \quad
(x_1, \cdots, x_k), \, (y_1, \cdots, y_k) \in \R^k.
\end{equation*}

\begin{prop} \label{proptau}
Let the conditions of Theorem \ref{Xexists}(i) be satisfied, $p\ge 1,$ and let
$\{X_t\}, \{r_t\}$ be stationary $L^p$-solutions of \eqref{Xform},  \eqref{genformQ},  respectively. In addition, assume that $b_j$ satisfy \eqref{bcond} with $\gamma >1$.  Then
\begin{equation}\label{taubdd}
\tau_p(n; \{X_j\}) = O(n^{-\gamma +1}), \qquad \tau_p(n; \{r_j\}) = O(n^{-\gamma +1}).
\end{equation}
\end{prop}
(\cite{ded2007}, Thm 1) together with Proposition \ref{proptau} imply the following
CLT for the empirical distribution functions
$F^X_n(u) := n^{-1}  \sum_{t=1}^n \1 (X_t \le u), \, F^r_n(u) $  $ := n^{-1}  \sum_{t=1}^n \1 (r_t \le u), \, u \in \R $ of stationary
solutions $\{X_t\}, \{r_t\}$  of \eqref{Xform},  \eqref{genformQ}.  Let $F^X (u) := \P(X_0 \le  u), F^r(u) :=  \P (r_0 \le u)$ be
the corresponding distribution functions. See \cite{ded2007} for the definition of weak convergence in the space $\ell^\infty (\R)$ of
all bounded functions on $\R$.

\begin{cor}  Let the conditions of Proposition \ref{proptau} hold with $p=1$ and $\gamma > 5$. Moreover, assume
that $F^X, F^r $ have bounded densities. Then
$\{ \sqrt{n} (F^X_n(u) - F^X(u)), u \in \R \} $ and $\{ \sqrt{n} (F^r_n(u) - F^r(u)), u \in \R \} $ converge
weakly in $\ell^\infty(\R)$ as $n \to \infty $ towards Gaussian processes on $\R$ with zero mean and respective
covariance functions $$\sum_{k\in \Z} {\rm Cov}(\1(  X_0 \le u), \1(X_k \le u)),\ \mbox{ and }\
\sum_{k\in \Z} {\rm Cov}(\1(r_0 \le u), \1(r_k \le u)) .$$
\end{cor}
\begin{rem}  {\rm
In \cite{ded2007} a general tightness condition is proposed in Proposition 6 for alternative classes of functions besides indicators of half-lines. Conditions are not immediate to check which explains that we restricted to the case of empirical cdf.}
\end{rem}

\section{Strong dependence}

The term {\it strong dependence} or {\it long memory} usually refers to stationary process
$\{y_t, t\in \Z \} $ whose covariance decays slowly with the lag so that its absolute series diverges:
$\sum_{k=1}^\infty $  $ |{\rm Cov}(y_0,y_k)| = \infty$. Since the variance of $\sum_{k=1}^n y_t$ usually grows faster than $n$ under
long memory, Donsker's invariance principle in \eqref{clt}  is no more valid and the limit of the partial  sums process, if exists,
might be quite complicated. Probably, the most important model of long memory processes
is the linear, or moving average, process $y_t = \sum_{s\le t}  b_{t-s} \zeta_s$, where $\{\zeta_s\}$ is
an i.i.d. sequence with zero  mean and finite variance, and the moving-average coefficients $b_j$ decay as in
\eqref{Qbc} below. Various generalizations of the linear model were studied in \cite{bai2008}, \cite{gru2014} and other works.
See the monograph \cite{gir2012} for a discussion and applications of long memory processes.

It is natural to expect that the `long memory' asymptotics of $b_j$ in \eqref{Qbc} induces some kind of long memory
of solutions $\{r_t\}, \{X_t\} $ of \eqref{genform}, \eqref{Xform}, under general assumptions
on $Q$. Concerning the latter process, this is true indeed as shown in the following theorem.

\begin{theo} \label{propXlong} Let $\{X_t\}$ be a stationary $L^2$-solution of \eqref{Xform}, where
\begin{equation}\label{Qbc}
b_j \ \sim \ \beta \, j^{d-1} \qquad (\exists \ 0< d < 1/2, \ \beta >0)
\end{equation}
and $Q$ satisfies the Lipschitz condition in \eqref{QLip} with $\lip^2_Q B^2 = \lip^2_Q \sum_{j=1}^\infty b^2_j < 1 $. Then
\begin{eqnarray}\label{Xcov}
&&{\rm Cov}(X_0, X_t)\ \sim\  \lambda^2_1 t^{2d-1}, \quad t \to \infty \qquad \text{and} \\
&&n^{-d - (1/2)} \sum_{t=1}^{[n\tau]} X_t \ \to_{D[0,1]} \ \lambda_2 B_{d + (1/2)}(\tau), \nn
\end{eqnarray}
where $B_{d+(1/2)} $ is a fractional Brownian motion with ${\rm Var}(B_{d+(1/2)}(\tau)) = \tau^{2d+1}$ and
$\lambda_1^2 := \beta^2 B(d,1-2d) \E Q^2(a+  X_0), \, \lambda^2_2 := \lambda_1^2/d(1+2d) $.
\end{theo}

Clearly, properties as in \eqref{Xcov} do  not hold for $\{r_t = \zeta_t Q(a+ X_t)\} $ which is an uncorrelated
martingale difference sequence. Here, long memory should appear in the behavior of the volatility $\sigma_t = Q(a + X_t), $
being `hidden' inside of nonlinear kernel $Q$. The last fact makes it much harder to prove it rigorously.
In the rest of the paper we  restrict ourselves to the quadratic model with $Q^2(x) = c^2 + x^2, $ or
\begin{equation}  \label{rsqr}
r_t \ = \  \zeta_t \sqrt{c^2 + \Big(a + \sum_{s<t} b_{t-s} r_s\Big)^2}, \qquad t \in \Z
\end{equation}
as in  \eqref{newQ}, where   (recall) $\{\zeta_t \}$ are standardized i.i.d. r.v.s, with zero mean and unit variance,
and $b_j, j\ge 1 $ are real numbers satisfying \eqref{Qbc}.

The following theorem  shows that under some additional conditions
the squared process  $\{r^2_t \}$ of \eqref{rsqr}
has similar long memory properties as $\{X_t\}$ in
Theorem \ref{propXlong}.
For the LARCH model (see Example 1 above) similar results were obtained in
(\cite{gir2000}, Thm. 2.2, 2.3).

\begin{theo} \label{long} Let $\{r_t\} $ be a stationary $L^2$-solution of \eqref{rsqr} with $b_j$ satisfying \eqref{Qbc}
and $B^2 = \sum_{j=1}^\infty b^2_j < 1$. Assume in addition that
$\mu_4 =\E [\zeta^4_0] < \infty $ and $\E r^4_t < \infty $.
Then 
\begin{equation} \label{cov2r}
{\rm Cov}(r^2_0, r^2_t) \  \sim \  \kappa^2_1 t^{2d-1 }, \qquad  t \to \infty
\end{equation}
where $\kappa^2_1 :=  \big(\frac{2a \beta}{1-B^2} \big)^2 B(d, 1-2d) \E r^2_0 $.   Moreover,
\begin{equation}\label{lim2r}
n^{-d-1/2} \sum_{t=1}^{[n\tau]} (r^2_t - \E r^2_t) \ \to_{D[0,1]} \  \kappa_2 B_{d +  (1/2})(\tau),
 \qquad  n \to \infty,
\end{equation}
where $B_{d +  (1/2)}$ is a fractional Brownian motion as in \eqref{Xcov}
and $\kappa^2_2 := \kappa^2_1/ (d(1+2d))$.

\end{theo}

\begin{rem} {\rm  Recently, Grublyt\.e and \v Skarnulis \cite{gru2015a} extended Theorem \ref{long} and some
other results of this paper to a more general QARCH model:
\begin{equation}  \label{qgarch}
r_t \ = \  \zeta_t \sigma_t,  \qquad \sigma^2_t \ = \ c^2 + \big(a + \sum_{s<t} b_{t-s} r_s\big)^2 + \gamma \sigma^2_{t-1}, \qquad t \in \Z
\end{equation}
involving lagged variable $\sigma^2_{t-1}$. 
For parametric `long memory' coefficients $b_j = \beta j^{d-1} $ as in \eqref{Qbc}
Grublyt\.e et al. \cite{gru2015b} discussed quasi-maximum
likelihood (QML) estimation of parameters $a, \beta, c, d, \gamma $ from observed sample  $r_t, 1\le t \le n$ satisfying the equations in \eqref{qgarch},
and proved consistency and asymptotic normality of the QML estimates. Related results for  parametric  LARCH
model with long memory were obtained in    Beran and  Sch\"utzner \cite{bera2009}.

}
\end{rem}

\section{Leverage}

Given a stationary conditionally heteroscedastic time series $\{r_t \}$ with $\E |r_t|^3 < \infty $ and
conditional variance $\sigma^2_t = {\rm Var}(r^2_t\,|\;r_s, s<t)$, leverage (a tendency of $\sigma^2_t$ to move
into the opposite direction as $r_s$ for $s<t$) is usually measured
by the covariance $h_{t-s} =  {\rm Cov}(\sigma^2_t, r_s). $ Following
\cite{gir2004}, we say that {\it $\{r_t\} $ has leverage of order $k$} $(1\le k <\infty)$ (denoted by
$\{r_t\} \in \ell(k)$)
whenever
\begin{equation}
h_j <0, \qquad 1\le j \le k.
\end{equation}
Note for $\{r_t\}$ in \eqref{genform},
\begin{equation}\label{hfunc}
h_j = \E [r^2_j r_0], \qquad j=0,1, \dots
\end{equation}
is the mixed moment function. Below, we show that
in the case of the quadratic
$\sigma^2_t$ in \eqref{newQ}, viz.,
\begin{equation}\label{root}
r_t = \zeta_t \sigma_t, \qquad \sigma_t = \big(c^2 + \big(a + \sum_{j=1}^\infty b_j r_{t-j}\big)^2 \big)^{1/2}
\end{equation}
and $\mu_3 = \E [\zeta^3_0] = 0$, the function
$h_j$ in \eqref{hfunc} satisfies a linear equation in \eqref{heq1}, below,
which can be analyzed and
the leverage effect for $\{r_t\}$ in \eqref{root}
established in spirit of  \cite{gir2004}. \\
Let $L^2(\Z_+) $  be the Hilbert space of all real sequences $\psi = (\psi_j, j \in \Z_+), \Z_+ :=
\{1,2,\cdots \}$ with finite norm $\|\psi\| := (\sum_{j=1}^\infty \psi_j^2)^{1/2} < \infty. $
As in the previous sections, let $B := \big(\sum_{j=1}^\infty b^2_j\big)^{1/2}$
and assume that $\{\zeta_t\}$ is an i.i.d. sequence
with zero mean and unit variance;  $\mu_i := \E \zeta_0^i, \, i=1,2,\cdots $.  \\
The following
theorem  establishes a criterion for the presence or absence of leverage in model \eqref{root},
analogous to the Thm. 2.4 in \cite{gir2004}.
\\
We also note that the proof of
Theorem \ref{lev} is  simpler than that of the above mentioned theorem, partly because
of the assumption $\mu_3 = 0$ used in the derivation of equation  \eqref{heq1}.

\begin{theo}\label{lev}
Let $\{r_t\} $ be a stationary $L^2$-solution of \eqref{rsqr} with $\E |r_0|^3 < \infty, |\mu|_3 < \infty$.  Assume in addition that
$B^2<1/5 $ and  
$\mu_3 = \E \zeta^3_0 = 0. $
Then for any fixed $k$ such that $1\le k \le \infty$:

\smallskip

\noi (i) if $a b_1<0$, $ab_j\le 0, \, j=2, \cdots, k$, then $\{r_t\}\in \ell(k)$

\smallskip

\noi (ii) if $a b_1>0$, $ab_j\ge 0, \, j=2, \cdots, k$, then $h_j>0, \, j=1, \cdots, k.$
\end{theo}

Particularly, for the Asymmetric ARCH(1) in \eqref{M11} with $\E |r_0|^3 < \infty,
\E \zeta_0^3 = 0$ the leverage function is $h_j = 2m_2 a b^{2j-1} $, see \eqref{m3},
and $\{r_t\} \in \ell(k)$ is equivalent to $ab < 0$.
Apparently, condition $B^2 < 1/5$ is not necessary for the statement of Theorem \ref{lev} although
a similar condition appears in the study of the leverage effect in  (\cite{gir2004}, (51)).

\section{Proofs}

\noi {\it Proof of Proposition \ref{Xreq}.}
(i) Since $\{X_t\}$ is predictable and $Q$ satisfies \eqref{Qnel} so
\begin{eqnarray*}
\E |r_t|^p&=&|\mu|_p
\E |Q(a+ X_t)|^p \\
&\le&|\mu|_p \E |c_1^2 + c_2^2 (a + X_t)^2|^{p/2} \\
&\le&C(1+ \E |X_t|^p) < C < \infty,
\end{eqnarray*}
proving \eqref{rX}.
Moreover, if $p> 1$ then
$\E [r_t |{\cal F}_{t-1}] = 0$ is a stationary martingale difference sequence.
Hence by Proposition \ref{Yp}, the series in \eqref{Xdef} converges in $L^p$ and satisfies
\begin{eqnarray*}
\E |X_t|^p &\le&C\left\{\begin{array}{ll}
\sum_{j=1}^\infty |b_j|^p, &0< p \le 2 \\
\big(\sum_{j=1}^\infty b_j^2\big)^{p/2}, &p > 2
\end{array}
\right\} \ = \ CB_p \ < \ \infty.
\end{eqnarray*}
In particular, $\zeta_t Q(a + \sum_{s<t} b_{t-s} r_s) = \zeta_t Q(a+ X_t) = r_t $
by the definition of $r_t$.  Hence, $\{r_t \}$ is a $L^p$-solution of (\ref{genformQ}).
Stationarity of   $\{r_t \}$ follows
from stationarity of $\{X_t \}$.  \\
Relations \eqref{genformvol} follow from $\E [\zeta_t| {\cal F}_{t-1}] = 0, \
\E [|\zeta_t|^p| {\cal F}_{t-1}] = |\mu|_p, \, p>1$, and the fact that $X_t$ is
${\cal F}_{t-1}$-measurable.

\smallskip

\noi (ii) Since $\{r_t\}$ is a $L^p$-solution  of \eqref{genformQ}, so $r_t = \zeta_t Q(a + X_t)$ with $X_t$ defined
in \eqref{Xdef}, and $\{X_t\}$ satisfy \eqref{Xform}, where the series converges in $L^p$. Also note that  $\{X_t \}$ is predictable.
Hence,  $\{X_t \}$ is a $L^p$-solution  of \eqref{Xform}.
By \eqref{Qnel},
\begin{eqnarray*}
\E|r_t|^p&=&|\mu|_p \E |Q(a+ X_t)|^p  \ \le\  |\mu|_p \E |c_1^2 + c_2^2 (a + X_t)^2|^{p/2} \
\ \le \  C(1+ \E |X_t|^p) \ < \  C.
\end{eqnarray*}
It also readily follows that, for $p> 1$,
$\{r_t, {\cal F}_t, t \in \Z\}$ is a  martingale difference sequence.
Hence, by the moment inequality in
\eqref{rosen},
\begin{eqnarray}\label{rosen1}
\E |X_t|^p &\le&K_p\left\{\begin{array}{ll}
\sum_{j=1}^\infty |b_j|^p \E|r_{t-j}|^p, &0< p \le 2 \\
\big(\sum_{j=1}^\infty b_j^2 \E^{2/p} |r_{t-j}|^p\big)^{p/2}, &p > 2
\end{array}
\right\} \ = \ CB_p \E |r_t|^p,
\end{eqnarray}
proving \eqref{Xr}. Stationarity of
 $\{X_t \}$ and \eqref{covX} are easy consequences of the above facts and stationarity of $\{r_t\}$.
 \hfill $\Box$

\vskip.2cm

\noi {\it Proof of Theorem \ref{Xexists}.} (i) For $n\in \N $ define a solution 
of  \eqref{Xform} with zero initial condition at $t \le - n $ as
\begin{eqnarray}
X^{(n)}_t&:=&\begin{cases}
0, &t \le -n, \\
\sum_{s= -n}^{t-1}b_{t-s} \zeta_s Q(a + X^{(n)}_s), &t> -n, \quad t \in \Z.
\end{cases}\label{Xndef}
\end{eqnarray}
Let us show that $\{X^{(n)}_t \}$ converges in $L^p$ to a stationary $L^p$-solution $\{X_t\}$
as $n \to \infty $. First, let $0< p \le 2$.
Let $m > n \ge 0$.  Then by inequality \eqref{rosen} for any $ t>  -m$ we have that
\begin{eqnarray*}
\E |X^{(m)}_t - X^{(n)}_t|^p&=&K_p|\mu|_p\Big\{\sum_{-m \le s<-n}  |b_{t-s}|^p
\E |Q(a + X^{(m)}_s)|^p\\
&+&\sum_{-n \le s < t} |b_{t-s}|^p  \E \big|Q(a + X^{(n)}_s) - Q(a + X^{(m)}_s)\big|^p\Big\}  \\
&=:& K_p|\mu|_p\{S'_{m,n} +  S''_{m,n}\}.
\end{eqnarray*}
Let $\chi_p(n) := \sum_{j=n}^\infty |b_j|^p$. \\
From the bound $|a + x|^2 \le (2a^2/\epsilon) + (1+\epsilon)x^2,$ valid for any $0< \epsilon < 1/2, $ $x\in \R$ and $a\ge 0,$ it follows that
\begin{eqnarray*}
\big|c_1^2 + c_2^2 (a + X^{(m)}_s)^2)\big|^{p/2}
&\le&c_1^p + c_2^p |(a + X^{(m)}_s)^2|^{p/2} \\
&\le&C(c_1,c_2) + c_2^p \,(1+ \epsilon)^{p/2}|X^{(m)}_s|^{p} \\
&\le&C(c_1,c_2) + c_3^p |X^{(m)}_s|^{p},
\end{eqnarray*}
with $c_3 > c_2 > \lip_Q$ arbitrarily close to $\lip_Q$. Then using (\ref{Qnel}) we obtain
\begin{eqnarray*}
S'_{m,n}
&\le&\sum_{-m\le s <-n} |b_{t-s}|^p
\E \big|c_1^2 + c_2^2 (a + X^{(m)}_s)^2\big|^p\\
&\le&C(Q)K_p |\mu|_p \chi_p(t+n) + c^p_3 \sum_{-m \le  s < -n}  |b_{t-s}|^p \E |X^{(m)}_s - X^{(n)}_s|^p, \\
S''_{m,n}
&\le&\lip^p_Q\sum_{-n \le s < t} |b_{t-s}|^p  \E \big|X^{(n)}_s - X^{(m)}_s\big|^p.
\end{eqnarray*}
Consequently,
\begin{eqnarray*}
\E |X^{(m)}_t - X^{(n)}_t|^p&\le&C(Q)K_p |\mu|_p\chi_p(t+n) +  K_p|\mu|_p c^p_3 \sum_{-m \le s < t} |b_{t-s}|^p \E \big|X^{(n)}_s - X^{(m)}_s\big|^p.
\end{eqnarray*}
Iterating the above inequality, we obtain
\begin{eqnarray}\label{series}
\E |X^{(m)}_t - X^{(n)}_t|^p&\le& C(Q)K_p|\mu|_p \Big\{\chi_p(t+n) + \sum_{k=1}^\infty  (K_p |\mu|_k c^p_3)^k \\
&&\times \sum_{-m \le  s_k < \cdots < s_1 < t} |b_{t-s_1}|^p  |b_{s_1-s_2}|^p \cdots |b_{s_{k-1}-s_k}|^p  \chi_p(s_k+n)\Big\}.\nn
 \end{eqnarray}
Since $K_p |\mu|_p c^p_3 B_p < 1 $ by \eqref{cQB}
and $\sup_{s\ge 1} \chi_p(s) \le B_p<\infty$, the  series on the r.h.s. of
\eqref{series}
is bounded uniformly in
$m, n $ and tends to  zero as $m, n \to \infty$ by the dominated convergence theorem.
Hence, there exist $X_t, t \in \Z$ such that
\begin{equation}\label{Xnconv}
\lim_{n\to  \infty} \E |X_t  - X^{(n)}_t|^p \ = \  0, \qquad \forall \ t \in \Z.
\end{equation}
Note  that $\{X_t\}$ is predictable and 
\begin{eqnarray*}
\E |X_t|^p\ =\ \lim_{n \to \infty}  \E |X^{(n)}_t|^p
&\le&\frac{C(Q) K_p |\mu|_p B_p }{1 - K_p |\mu|_p c^{p}_3 B_p} \
\le\ \frac{C(p,Q) |\mu|_p B_p }{1 - K_p |\mu|_p \lip^{p}_Q B_p},
\end{eqnarray*}
where the last inequality follows by taking $c_3 >\lip_Q$ sufficiently close to $\lip_Q$. \\
We also have by \eqref{X2mom} and \eqref{QLip} that
\begin{eqnarray*}
&&\E \big|\sum_{s< t} b_{t-s} \zeta_s Q(a + X_s) - \sum_{s= -n}^{t-1}b_{t-s} \zeta_s Q(a + X^{(n)}_s)\big|^p \\
&&=\  \E \big|\sum_{s< -n} b_{t-s} \zeta_s Q(a + X_s) + \sum_{s= -n}^{t-1}b_{t-s} \zeta_s
(Q(a+X_s) - Q(a + X^{(n)}_s))\big|^p\\
&&\le\  K_p |\mu|_p\Big\{\sum_{ s < -n} |b_{t-s}|^p \E \big|Q(a + X_s)\big|^p +
\sum_{-n \le s <t} |b_{t-s}|^p \E \big|Q(a + X_s) -  Q(a + X^{(n)}_s)\big|^p \Big\} \\
&&\le\  C\Big( \sum_{ s < -n} |b_{t-s}|^p +\sum_{s <t} |b_{t-s}|^p \E \big|X_s -  X^{(n)}_s\big|^p \Big) \ \to \  0
\end{eqnarray*}
as  $n \to \infty$.
Whence and from \eqref{Xndef} it follows that $\{X_t \} $ is a stationary $L^p$-solution
of  \eqref{Xform} satisfying  \eqref{X2mom}. \\
To show the uniqueness of stationary $L^p$-solution of \eqref{Xform},
let $\{ X'_{t}\}, \{X''_{t}\}$ be  two such solutions of \eqref{Xform}, and
$m_p(t) := \E |X'_t - X''_t|^p. $  \\
Then $\sup_{t\in \Z} m_p(t) \le  M < \infty $ and
$ m_p(t) \le K_p |\mu|_p \lip^p_Q \sum_{s<t} |b_{t-s}|^p m_p(s)$ follows by \eqref{QLip}.  Iterating the last equation we obtain that
$m_p(t) \le  (K_p |\mu|_p \lip^p_Q B_p)^k  M $ holds for all  $k \ge 1$, where $K_p |\mu|_p \lip^p_Q B_p < 1$. \\
Hence, $m_p(t) = 0$.  This proves part (i) for $0< p \le 2$. \\
The proof of  part (i) for $p> 2$ is analogous. Particularly, using \eqref{rosen} as in \eqref{rosen1},
we obtain
\begin{eqnarray*}
\E |X_t|^p &\le&K_p |\mu|_p\big(\sum_{s<t} b_{t-s}^2 \E^{2/p} |Q(a + X_s)|^p\big)^{p/2} \nn \\
&\le&K_p |\mu|_p \big(\sum_{s<t} b_{t-s}^2 (C(Q) + c_3^p \E |X_s|^p )^{2/p} \big)^{p/2} \nn \\
&\le&K_p|\mu|_p B_p (C(p,Q) + c^p_3 \sup_{s\in \Z} \E|X_s|^p) \label{naujas}
\end{eqnarray*}
implying $(1- K_p |\mu_p| c_3^p B_p )\sup_{t\in \Z} \E |X_t|^p \le C(p,Q)  |\mu|_p B_p $ and
hence the bound in \eqref{X2mom} for $p> 2$, by taking $c_3 $ sufficiently close to $\lip_Q$.
This proves part (i).

\smallskip

\noi (ii) Note that $Q(x) =  \sqrt{c_1^2 + c_2^2 x^2}$ is a Lipschitz function and satisfies \eqref{QLip} with $\lip_Q = c_2$. Hence by $K_2 =1 $ and part (i), a unique
$L^2$-solution $\{X_t\}$ of \eqref{Xform} under the condition $c^2_2 B_2  < 1$ exists. To show
the necessity of the last condition,
let $\{X_t \} $ be a stationary $L^2$-solution of  \eqref{Xform}.  Then
\begin{eqnarray*}
\E X^2_t&=&\sum_{s<t}  b^2_{t-s} \E Q^2(a + X_s)\\
&=&\sum_{s<t}  b^2_{t-s} \E \big(c_1^2 + c_2^2 (a+ X_s)^2 \big) \\
&=&B_2\big(c_1^2 + c^2_2 (a^2 + \E X^2_t)\big) \ > c^2_2 B_2 \E X^2_t
\end{eqnarray*}
since $a \neq 0$.  Hence, $c^2_2 B_2  < 1$ unless $\E X^2_t = 0 $,
or  $\{ X_t = 0\}$ is a trivial process.
Clearly, \eqref{Xform} admits a trivial solution if and only if $0= Q(a) = \sqrt{c_1^2 + c_2^2 a^2} = 0, $ or
$c_1 = c_2 = 0$. This proves part (ii) and the theorem. \hfill  $\Box$

\vskip.2cm

\noi The proofs of Proposition \ref{propomega} and Theorem \ref{long} use the  following lemma.

\vskip.2cm

\begin{lem} \label{lemmaA} For $\alpha_j \ge 0, j=1,2, \cdots, $ denote
\begin{eqnarray} \label{Ak}
A_k \ := \  \alpha_k + \sum_{0< p <k} \sum_{0< i_1 <  \dots < i_p <  k} \alpha_{i_1} \alpha_{i_2 -i_1}  \cdots \alpha_{i_p - i_{p-1}} \alpha_{k-i_p},
\qquad k =  1,2, \cdots.
\end{eqnarray}
Assume that $A :=  \sum_{j=1}^\infty \alpha_j < 1 $ and
\begin{equation}\label{alpha}
\alpha_j \ \le  \  c \, j^{-\gamma},  \qquad (\exists  \  c >0, \ \gamma >  1).
\end{equation}
Then there exists  $C >0$ such  that for any $k \ge 1 $
\begin{equation} \label{Ak1}
A_k  \  \le \   C k^{-\gamma}.
\end{equation}
\end{lem}

\noi {\it  Proof.} We have $A_k = \sum_{0\le p <k} A_{k,p}$, where
$$
A_{k,p} := \sum_{0< i_1 <  \dots < i_p <  k} \alpha_{i_1} \alpha_{i_2 -i_1}  \cdots \alpha_{i_p - i_{p-1}} \alpha_{k-i_p} \quad (p\ge 1),
\quad A_{k,0} :=  \alpha_k
$$
is the inner sum in \eqref{Ak}. W.l.g., assume $c \ge 1 $ in \eqref{alpha}.
Let us prove that  there exists $\lambda >0$ such that
\begin{equation}\label{Akp}
A_{k,p}  \  \le \   c (p+2)^\lambda A^{p+1}  k^{-\gamma}, \qquad \forall \   0 \le p < k < \infty.
\end{equation}
Since $A<  1$,  so \eqref{Akp} and $\sum_{p>0}  (p+2)^\lambda A^{p+1}  < \infty $ together
imply \eqref{Ak1}.

By dividing both  sides of \eqref{Akp} by $A^{p+1}$, it suffices to show \eqref{Akp} for  $A =1 $.
The proof  uses induction on $p $.  Clearly, \eqref{Akp} holds for $p=0$. To prove the induction
step $p-1 \to p \ge 1$, note
\begin{eqnarray} \label{Ak2}
A_{k,p}&=&\sum_{0< i <k} \alpha_i A_{k-i,p-1}  = \sum_{\frac{k}{p+1}< i < k}  \alpha_i A_{k-i,p-1} +
\sum_{k - \frac{k}{p+1} \le   k-i < k}  \alpha_i A_{k-i,p-1}.
\end{eqnarray}
Here,  $\alpha_i \1( i> \frac{k}{p+1}) \le c i^{-\gamma} \1( i> \frac{k}{p+1}) \le c (p+1)^\gamma k^{-\gamma} $ and,
similarly,  by the inductive assumption
$$
A_{k-i,p-1} \1( k-i \ge k - \frac{k}{p+1})  \le  c (p+1)^\lambda (k - \frac{k}{p+1})^{-\gamma} =
c (p+1)^\lambda \big(\frac{p+1}{p}\big)^\gamma k^{-\gamma}.
$$
Assumption $A=1$ implies $\sum_{k>0}  A_{k,p} = 1$ for any $p\ge  0$. Using the above facts
from \eqref{Ak2}  we obtain
\begin{eqnarray*}
A_{k,p}&=&\frac{c (p+1)^\gamma}{k^\gamma} \sum_{k/(p+1)< i < k}  A_{k-i,p-1} +
\frac{c (p+1)^\lambda}{k^\gamma} \big(\frac{p+1}{p}\big)^\gamma
\sum_{k - k/(p+1) \le   k-i < k}  \alpha_i  \\
&\le&c\big((p+1)^\gamma +  (p+1)^\lambda \big(\frac{p+1}{p}\big)^\gamma\big)k^{-\gamma}.
\end{eqnarray*}
Hence the proof of the induction step $p-1 \to p \ge 1$ amounts to verifying the inequality
$(p+1)^\gamma +  (p+1)^\lambda \big(\frac{p+1}{p}\big)^\gamma \le (p+2)^\lambda $, or
\begin{equation} \label{induc}
n^\gamma +  n^\lambda \big(\frac{n}{n-1}\big)^\gamma \ \le \ (n+1)^\lambda, \qquad n=2,3, \dots.
\end{equation}
The above inequality holds with
$\lambda = 3 \gamma $. Indeed,
\begin{eqnarray*}
n^\gamma +  n^{\lambda} \big(\frac{n}{n-1}\big)^\gamma &=&
n^{\lambda} (n^{-2 \gamma } +  \big(\frac{n}{n-1}\big)^\gamma )\le
 n^{\lambda} (n^{-2} +  \big(\frac{n}{n-1}\big))^\gamma \\
&\le&  n^{\lambda} (1+  \frac{1}{n-1}+ \frac{1}{n^{2}} )^\gamma \le  n^{\lambda} (1+  \frac{3}{n}+ \frac{3}{n^{2}} +\frac{1}{n^3})^\gamma
= (n+1 )^{\lambda},
\end{eqnarray*}
proving \eqref{induc} and the lemma, too. \hfill $\Box$

\bigskip

\noi {\it Proof of Proposition \ref{propomega}.} We will give the proof for $p \ge 2 $ only as the proof for
$p\in [1,2]$ is similar. \\
Following the notation in \eqref{omega2}, let $\{X'_t\}, \{r'_t\}$ be the corresponding processes (Bernoulli shifts)
of the i.i.d. sequence $\xi' := (\cdots, \zeta_{-1}, \zeta'_0, \zeta_1, \zeta_2, \cdots)$ with $\zeta_0$ replaced by its
independent copy $\zeta_0'$. Note that $X'_t = X_t \, (t \le 0), \, r'_t = r_t \, (t<0)$. We have
$\delta^2_2(k; \{X_t\} ) = (\E |X_k - X'_k|^p )^{2/p} = \|X_k - X'_k\|_p^2 $, where
$$
X_k - X'_k \ = \  \sum_{0<s<k} b_{t-s} (r_s - r'_s) + b_k (\zeta_0 - \zeta'_0) Q(a + X_0).
$$
Then with $v^2_p  := \| Q(a + X_0)\|^2_p$ using Rosenthal's inequality \eqref{rosen} similarly as in
the proof of Theorem \ref{Xexists} we obtain
\begin{eqnarray*}
\|X_k - X'_k\|_p^2&\le&K_p^{2/p}
\Big(\sum_{0<s<k} b^2_{k-s} \|r_s - r'_s\|_p^2 + \|\zeta_0 - \zeta'_0\|^2_p b^2_k v^2_p \Big)\\
&\le&K_p^{2/p}\Big(\sum_{0<s<k} b^2_{k-s} |\mu|_p^{2/p} \|Q(a+ X_s)- Q(a+ X'_s)\|_p^2 + 4 |\mu|_p^{2/p} b^2_k v^2_p \Big) \\
&\le&K_p^{2/p}  |\mu|^{2/p}_p \Big( \lip^2_Q \sum_{0<s<k} b^2_{k-s} \|X_s- X'_s\|_p^2 + 4 b^2_k v^2_p\Big).
\end{eqnarray*}
Let $\alpha_k := K_p^{2/p}  |\mu|^{2/p}_p \lip^2_Q b^2_k$. Iterating the last inequality we obtain
\begin{eqnarray*}
\delta^2_2(k; \{X_t\})&\le&\frac{4 v^2_p}{\lip^2_Q}  \big(\alpha_k + \sum_{0< s <k} \alpha_s \alpha_{k-s} + \cdots \big)\
=\ \frac{4 v^2_p}{\lip^2_Q} \cdot A_k,
\end{eqnarray*}
where $A_k$ is as in \eqref{Ak}. Since $A = \sum_{k>0} \alpha_k =  (K_p  \mu_p \lip^p_Q B_p)^{2/p} < 1 $ and
$\alpha_k \le C k^{-2\gamma}$, by Lemma \ref{lemmaA} we obtain
$ \delta_2(k; \{X_t\}) \le C k^{-\gamma}$,  proving
the first inequality in \eqref{deltas}. The proof of the second inequality in \eqref{deltas}
follows similarly using  $ \delta^2_p(k; \{r_t\}) = \|r_k - r'_k\|^2_p \le \lip^2_Q \mu_p^{2/p} \| X_k - X'_k\|^2_p
= \lip^2_Q \mu_p^{2/p} \delta^2_2(k; \{X_t\})$. \hfill $\Box$

\bigskip

\noi {\it Proof of Proposition \ref{proptau}.} We use the coupling inequality of \cite{ded2005} in \eqref{coup}, below, providing
a simple upper bound for $\tau$-coefficients. Let $\{y^*_j \} $ be distributed as $\{y_t\}$ and independent of $y_s, s \le 0$. Then
\begin{eqnarray}\label{coup}
\tau_p(\{y_{j_i}\}_{1\le i \le k})&\le&\sum_{i=1}^k \|y_{j_i} - y^*_{j_i}\|_p \quad
\text{and} \quad \tau_p(n,\{y_t\})\ \le \ \sup_{j \ge n}  \|y_j -  y_j^* \|_p.
\end{eqnarray}
To construct the coupling for $\{X_t\}$, let $\{ X^*_t\}$ be the corresponding
process (Bernoulli shift) of the i.i.d. sequence $\xi^* := (\cdots, \zeta^*_{-2}, \zeta^*_{-1}, \zeta_0, \zeta_1, \cdots ) $
with $(\zeta^*_s, s <0)$ an independent copy of $(\zeta_s, s<0)$.
Clearly, $\{X^*_t\}$ is distributed as $\{X_t\}$ and independent of $(X_s, s\le 0)$, the latter being measurable
w.r.t. $(\zeta_s, s < 0)$. Hence, the first relation in \eqref{taubdd}  follows from
\begin{equation}\label{XX}
\|X_n - X^*_n\|_p = O(n^{-\gamma+1}).
\end{equation}
Towards this end, consider `intermediate' i.i.d. sequence $\xi^*_i :=
(\cdots, \zeta^*_{-i-1}, \zeta^*_{-i}, \zeta_{-i+1}, \cdots, \zeta_0, $
$\zeta_1, \cdots), \, i\ge 1, \, \xi^*_1 := \xi^*$. Note
sequences $\xi^*_i$ and $\xi^*_{i+1}$ agree up to single entry. Let
$\{ X^*_{i,t}\}$ be the corresponding
Bernoulli shift of the i.i.d. sequence $\xi^*_i$.  By triangle inequality,
$\|X_n - X^*_n\|_p \le \sum_{i=1}^\infty \|X^*_{i,n} - X^*_{i+1,n}\|_p$. By stationarity and Proposition
\ref{propomega},
\begin{equation}\label{XXX}
\|X^*_{i,n} - X^*_{i+1,n}\|_p = \delta_p (n+i, \{X_t\}) \ \le \ C(n+i)^{-\gamma},
\end{equation}
where $\delta_p$ is defined in \eqref{omega2}. Clearly, \eqref{XXX} implies \eqref{XX}, proving
the first relation in  \eqref{taubdd}. Since $\tau_p (n, \{r_t\}) \le C_p \tau_p (n, \{X_t\})$, the second
relation in  \eqref{taubdd} follows. Proposition \ref{proptau} is proved. \hfill $\Box$

\bigskip

\noi{\it Proof of Theorem \ref{propXlong}.}  The first relation in \eqref{Xcov} follows from \eqref{covX} and
\eqref{Qbc}. The second relation in \eqref{Xcov} follows from a general result in Abadir et al. (\cite{aba2014}, Prop.3.1), using
the fact that $\{r_s\}$ in \eqref{Xform} is a stationary ergodic martingale difference sequence. \hfill $\Box$

\bigskip

\noi {\it Proof of  Theorem \ref{long}.} The proof of Theorem \ref{long} heavily relies on the decomposition
\begin{equation}\label{appX}
(r^2_t - \E r^2_t) - \sum_{s<t} b^2_{t-s} (r^2_s - \E r^2_s) \ = \ 2a X_t + Z_t,
\end{equation}
where $\{Z_t\}$ on the r.h.s. of \eqref{appX} is negligible so its memory intensity is less than the memory
intensity of the main term, $\{X_t\}$. Accordingly, $r^2_t - \E r^2_t = (1 - \sum_{j=1}^\infty b^2_j L^j )^{-1} \xi_t$
behaves like an AR($\infty$) process with long memory innovations $\xi_t := 2a X_t + Z_t \approx 2a X_t $. A rigorous
meaning to the above heuristic explanation is provided below. \\
By the definition of $r_t$ in \eqref{rsqr},
\begin{eqnarray}\label{Udef}
Z_t&:=&U_t + V_t, \qquad \text{where}  \nn \\
U_t&:=&(\zeta^2_t -1) Q^2(a + X_t), \nn \\
V_t&:=&X^2_t - \E X^2_t - \sum_{s<t} b^2_{t-s} (r^2_s - \E r^2_s) \nn \\
&=&2\sum_{s_2 < s_1  <t} b_{t-s_1} b_{t-s_2} r_{s_1} r_{s_2}. \label{Vdef}
\end{eqnarray}
Let us first check that the double series in \eqref{Vdef}
converges in mean square and \eqref{Vdef} holds.
Let
$$
X_{t,N} :=  \sum_{-N < s <t} b_{t-s} r_s, \qquad
V_{t,N} :=  2\sum_{-N< s_2 < s_1  <t} b_{t-s_1} b_{t-s_2} r_{s_1} r_{s_2},
$$
then $V_{t,N} = X^2_{t,N} - \E X^2_{t,N} - \sum_{-N< s<t} b^2_{t-s} (r^2_s - \E r^2_s)$ and, for $M > N$,
$$
\E (X^2_{t,N} - X^2_{t,M})^2 = \E (X_{t,N} - X_{t,M})^2(X_{t,N} + X_{t,M})^2
\le \|X_{t,N} - X_{t,M}\|^2_4 \|X_{t,N} + X_{t,M}\|^2_4.
$$
By Rosenthal's inequality in \eqref{rosen},
\begin{eqnarray*}
\|X_{t,N} + X_{t,M}\|^2_4&\le&C \sum_{-M < s < t}  b^2_{t-s} \ \le \ C \qquad  \text{and}  \\
\|X_{t,N} - X_{t,M}\|^2_4&\le&C \sum_{-M < s \le -N } b^2_{t-s} \ \to \ 0 \qquad  (N, M \to \infty).
\end{eqnarray*}
Therefore, $\lim_{N,M \to \infty} \E (X^2_{t,N} - X^2_{t,M})^2 = 0$. \\
The convergence of
$\E X^2_{t,N} $ and $\sum_{-N< s<t} b^2_{t-s} (r^2_s - \E r^2_s)$ in $L^2$  as $N \to \infty$ is easy. Hence,
$V_{t,N}, N \ge 1$ is a Cauchy sequence in $L^2$ and the double series in \eqref{Vdef} converges
as claimed above, proving  \eqref{Vdef}. \\
Let us prove that in the decomposition \eqref{appX},  $\{Z_t\} $
is negligible in the sense that its (cross)covariances decay faster as the covariance
of the main term, $\{X_t\}$, viz.,
\begin{equation}\label{Zrem}
\E [Z_t Z_0] =  o(t^{2d-1}), \qquad \E [X_t Z_0] =  o(t^{2d-1}), \qquad \E [Z_t X_0] =  o(t^{2d-1})
\end{equation}
as  $t \to \infty$. Note, for $t \ge 1$,
$\E [U_0 U_t] = \E [V_0 U_t] =0$ and $\E [V_t U_0]  = 2b_t \E [\zeta_0 (\zeta^2_0 -  1) Q^2(a+ X_0) \sum_{s_2 <0} b_{t-s_2} r_{s_2} ] =
O(b_t) = o(t^{2d-1})$.  Hence, the first relation in \eqref{Zrem} follows from
\begin{equation}\label{Vrem}
\E [V_t V_0] =  o(t^{2d-1}), \qquad t \to \infty,
\end{equation}
which is proved below. Since $\E [V^2_t] <  \infty, \E [V_t] = 0$ we can write the orthogonal
expansion
\begin{equation*}
V_t = \sum_{s < t} P_s V_t,
\end{equation*}
where $P_s V_t := \E [ V_t|{\cal  F}_s] - \E [V_t|{\cal F}_{s-1}]$ is the projection operator. \\
By orthogonality of $P_s$,
\begin{eqnarray*}
\big|\E V_0 V_t\big| &=&\big|\sum_{s<0} \E  [ (P_s V_0) (P_s V_t)] \big| \ \le \  \sum_{s<0} \|P_s V_0\|_2 \| P_s V_t \|_2.
\end{eqnarray*}
Relation \eqref{Vrem} follows from
\begin{eqnarray}\label{Pbdd}
 \| P_s V_0 \|^2_2 = o (b^2_{-s})  = o((-s)^{2(d-1)}), \qquad s \to  -\infty.
\end{eqnarray}
Indeed, if \eqref{Pbdd} is true then
\begin{eqnarray*}
\E V_0 V_t&=&o\big(\sum_{s<0} (-s)^{d-1} (t-s)^{d-1} \big) \ = \  o(t^{2d-1}), \qquad t \to \infty,
\end{eqnarray*}
proving \eqref{Vrem}. \\
Consider \eqref{Pbdd}.  We have by  \eqref{Vdef} and the martingale difference property of
$\{r_s\}$ that
\begin{eqnarray*}\label{Pform}
P_s V_0&=&2r_s b_{-s} \sum_{u < s} b_{-u} r_{u}
\end{eqnarray*}
and
\begin{eqnarray*}
\|P_s V_0\|_2^2&=&4 b^2_{-s} \E \big[r^2_s \big(\sum_{u < s} b_{-u} r_{u}\big)^2 \big] \
\le\  4b^2_{-s} \|r_s\|^2_4 \,  \big\| \sum_{u < s} b_{-u} r_{u} \big\|_4^2. \label{Pform2}
\end{eqnarray*}
By Rosenthal's inequality in \eqref{rosen},
\begin{eqnarray*}
\E \big| \sum_{u < s} b_{-u} r_{u} \big|^4&\le&C_4\big(\sum_{u< s} b^2_{-u}  (\E r^4_u)^{1/2} \big)^2 \ \le \
C  \big(\sum_{u> |s|}  u^{2(d-1)} \big)^2 \ = \ O(|s|^{2(2d-1)}) \ = \ o(1).
\end{eqnarray*}
Therefore,
\begin{eqnarray*}
\|P_s V_0\|_2^2
&\le&C|s|^{2(d-1) +  2d - 1} \ = \ o(|s|^{2(d-1)}), \label{Pform3}
\end{eqnarray*}
proving \eqref{Pbdd}, \eqref{Vrem}, and the first relation in \eqref{Zrem}.
The remaining
two relations  in \eqref{Zrem} follow easily, e.g.,
\begin{eqnarray*}
\E[ X_t Z_0]&=&b_t \E [r_0 (\zeta_0^2 -1) Q^2(a + X_0)] +
2\sum_{s_1 <0} b_{t-s_1} b_{-s_1} L_{s_1},
\end{eqnarray*}
where
\begin{eqnarray*}
L_{s_1}&:=&\E [r^2_{s_1} \sum_{s_2 < s_1} b_{-s_2} r_{s_2} ] \\
&\le&\E^{1/2} [r^4_{s_1}] \E^{1/2} \big[\big(\sum_{s_2<s_1} b_{-s_2} r_{s_2} \big)^2 \big] \nn \\
&=&O\big(\big(\sum_{s_2<s_1} b^2_{-s_2}\big)^{1/2}\big) \ = \ O( |s_1|^{d-(1/2)}), \qquad s_1 \to -\infty. \label{Lineq}
\end{eqnarray*}
Therefore
\begin{equation*}
\E[ X_t Z_0] \ = \ O(t^{d-1}) + \sum_{s_1 <0} (t-s_1)^{d-1} (-s_1)^{2d - (3/2)} \ = \  o(t^{2d-1}).
\end{equation*}
This proves \eqref{Zrem}. \\
Next, let us prove \eqref{cov2r}. Recall the decomposition  \eqref{appX}.
Denote $\xi_t := 2aX_t + Z_t $, then  \eqref{appX} can be rewritten
as $(r^2_t - \E r^2_t)- \sum_{s<t} b^2_{t-s} (r^2_s - \E r^2_s) = \xi_t  $, or
\begin{equation}\label{r22inv}
r^2_t - \E r^2_t \ = \ \sum_{i=0}^\infty \phi_i \xi_{t-i},  \qquad t \in \Z,
\end{equation}
where $\phi_j \ge 0, j\ge 0$ are the coefficients of the power series
$$
\Phi(z) := \sum_{j=0}^\infty \phi_j z^j
=  (1 - \sum_{j=1}^\infty b^2_j z^j)^{-1}, \quad  z \in \C, \ |z| < 1
$$
given by $\phi_0 := 1$,
\begin{eqnarray}\label{phi}
\phi_j&:=&b_j^2+ \sum_{0<k<j} \sum_{0<s_1 < \dots < s_k < j} b^2_{s_1} b^2_{s_2 -s_1} \cdots b^2_{s_k -s_{k-1}} b^2_{j-s_k}, \quad j\ge 1.
\end{eqnarray}
From \eqref{Qbc} and Lemma \ref{lemmaA} we infer that
\begin{equation}\label{philim}
\phi_t =  O(t^{2d-2}), \qquad t \to \infty,
\end{equation}
in particular, $\Phi(1) = \sum_{t=0}^\infty \phi_t = 1/(1- B^2) < \infty$ and the r.h.s. of \eqref{r22inv}
is well-defined. Relation \eqref{Zrem} implies that that
\begin{equation}\label{gammalim}
\gamma_t := {\rm Cov}(\xi_0, \xi_t) \sim  4a^2 {\rm Cov}(X_0, X_t) \sim 4a^2 \kappa^2_3 t^{2d-1}, \ t \to \infty
\end{equation}
with $\kappa^2_3 = \beta ^2 B(d, 1-2d)\E r^2_0$. Let us show that
\begin{eqnarray}\label{covvv}
{\rm Cov}(r^2_t, r^2_0)&=&\sum_{i,j=0}^\infty \phi_i \phi_j \gamma_{t-i+j}  \
\sim \ \Phi^2(1) \gamma_t, \qquad t \to \infty.
\end{eqnarray}
With \eqref{gammalim} in mind,  \eqref{covvv} is equivalent to
\begin{equation}\label{Jlim}
J_t := \sum_{i,j=0}^\infty \phi_i \phi_j (\gamma_{t-i+j} - \gamma_t) \ = \  o(t^{2d-1}).
\end{equation}
For a large $L >0$, split $J_t = J'_{t,L} + J''_{t,L}$, where
\begin{equation*}
J'_{t,L} := \sum_{i,j>0: |j-i| \le L} \phi_i \phi_j (\gamma_{t-i+j} - \gamma_t),  \qquad
J''_{t,L} := \sum_{i,j>0: |j-i| >L} \phi_i \phi_j (\gamma_{t-i+j} - \gamma_t).
\end{equation*}
Clearly, \eqref{Jlim} follows from
\begin{eqnarray}\label{JlimL}
t^{1-2d}J'_{t,L}&=&o(1) \quad \forall  \  L >0, \quad \text{and}  \quad
\lim_{L \to \infty} \limsup_{t\to \infty} t^{1-2d} J''_{t,L} \ = \ 0.
\end{eqnarray}
The first relation in \eqref{JlimL} is immediate from \eqref{gammalim} since the latter implies
$\gamma_{t+k} - \gamma_t = o(t^{2d-1})$ for any $k$ fixed.  \\
With \eqref{philim} and \eqref{gammalim} in mind,
the second relation in  \eqref{JlimL} follows from
\begin{eqnarray}\label{JlimL1}
\lim_{L \to \infty} \limsup_{t\to \infty} t^{1-2d} \bar J_{t,L}
 \ = \ 0.
\end{eqnarray}
where $\bar J_{t,L} := \sum_{i,j>0: |j-i| >L} i^{2d-2} j^{2d-2} (t^{2d-1}+|t+j-i|^{2d-1}_+)$ and where $k_+^{2d-1} := \min (1, k^{2d-1}),
k \in \Z_+$.
 \\
Split the last sum according to whether
$|t+j-i| \ge t/2$, or $|t+j-i| < t/2$. \\
Then
\begin{eqnarray*}
\bar J'_{t,L}&:=&\sum_{i,j>0: |j-i| >L, |t+j-i| \ge t/2} i^{2d-2} j^{2d-2} (t^{2d-1} +|t+j-i|^{2d-1} )\\
&\le&C t^{2d-1} \sum_{i,j >0:  |j-i| > L} i^{2d-2} j^{2d-2} \le C t^{2d-1} L^{2d-1}
\end{eqnarray*}
follows by  $ \sum_{i,j >0:  |j-i| > L} i^{2d-2} j^{2d-2} \le \sum_{0< i < L/2, j > L/2} i^{2d-2} j^{2d-2}
+ \sum_{i> L/2, j>0} i^{2d-2} j^{2d-2} = O(L^{2d-1})$. Therefore,
$\lim_{L \to \infty} \limsup_{t\to \infty} t^{1-2d} \bar J'_{t,L}
= 0. $ \\
Next, since  $|t+j-i| < t/2$ implies $i> t/2$, so with $k:= t+j-i$ we obtain
\begin{eqnarray*}
\bar J''_{t,L}
&\le&Ct^{2d-2} \sum_{i,j >0: |t+j-i| < t/2} j^{2d-2}  (t^{2d-1} +|t+j-i|_+^{2d-1}) \\
&\le&Ct ^{2d-2}\sum_{j>0} j^{2d-2} \sum_{|k| < t/2}  (t^{2d-1} +|k|_+^{2d-1}) \\
&\le&Ct ^{4d-2},
\end{eqnarray*}
implying $\limsup_{t\to \infty} t^{1-2d} \bar J''_{t,L} = 0$ for any $L>0$. This proves
 \eqref{JlimL},  \eqref{Jlim}, and  \eqref{covvv}. Clearly,
\eqref{cov2r} follows from  \eqref{covvv} and \eqref{gammalim}.

\smallskip

\noi It remains to show the invariance principle in \eqref{lim2r}. With \eqref{r22inv} in mind,
decompose $S_n(\tau):= \sum_{t=1}^{[n\tau]} (r^2_t - \E r^2_t)  = \sum_{i=1}^3 S_{ni}(\tau)$, where
\begin{eqnarray*}
S_{n1}(\tau)&:=&2a\Phi(1)\sum_{t=1}^{[n\tau]} X_t, \\
S_{n2}(\tau)&:=&\Phi(1) \sum_{t=1}^{[n\tau]}Z_t, \\
S_{n3}(\tau)&:=&\sum_{t=1}^{[n\tau]} \sum_{i=0}^\infty \phi_i(\xi_{t-i}-\xi_t).
\end{eqnarray*}
Here, $\E S^2_{n2}(\tau) = o(n^{2d+1})$ follows from  \eqref{Zrem}. Consider
\begin{eqnarray*}
\E S^2_{n3}(\tau)
&:=&\sum_{t,s=1}^{[n\tau]} \sum_{i,j=0}^\infty \phi_i \phi_j \E (\xi_{t-i}-\xi_t)(\xi_{s-j}-\xi_s)
\ = \ \sum_{t,s=1}^{[n\tau]} \rho_{t-s},
\end{eqnarray*}
where $\rho_t :=  \sum_{i,j=0}^\infty \phi_i \phi_j (\gamma_{t+j-i} - \gamma_{t+j} - \gamma_{t-i} + \gamma_t)
= o(t^{2d-1}) $ follows similarly to  \eqref{JlimL}. Hence,
$S_{ni}(\tau) = o_p(n^{-d-1/2}), \, i=2,3$. The convergence
$n^{-d-1/2}S_{n1}(\tau)  \to_{D[0,1]}   \kappa_2 B_{d +  (1/2})(\tau) $ follows from
Theorem \ref{propXlong}. \\
This completes the proof of Theorem \ref{long}. \hfill $\Box$

\bigskip

\noi {\it Proof of Theorem \ref{lev}.}  Let us first prove that $\|h\| < \infty $.
Note  that
\begin{equation} \label{3mixed}
\lim_{n\to \infty} \E \big(\sum_{-n<s<t} b_{t-s} r_s \big)^2 r_0 \ = \  \E \big(\sum_{-\infty <s<t} b_{t-s} r_s \big)^2 r_0,
\end{equation}
which follows from the definition of $L^3$-solution of \eqref{rsqr} and Remark \ref{Lpsol}.
Then using \eqref{3mixed},
$\E r_t = \E [r_t^3]
= \E [r_t r_s] = 0, s < t $ 
we obtain
\begin{eqnarray}
h_j&=&\lim_{n \to \infty} \E\Big[\Big(c^2+a^2+2a \sum_{-n<s<t}b_{t-s} r_s +
\sum_{-n<s<t} b^2_{t-s} r_s^2 \nn  \\
&&\qquad + \ 2 \sum_{-n<s_2 < s_1<t} b_{t-s_1} b_{t-s_2}
r_{s_1} r_{s_2}\Big)r_{t-j} \Big] \nn \\
&=&2am_2 b_j + \sum_{t-j<s<t} b_{t-s}^2 h_{j+ s-t} +
2 b_j \lim_{n \to \infty} \E R_n(t,j) \label{L1}
\end{eqnarray}
where  $R_n(t,j) :=  r^2_{t-j} \sum_{-n < s < t-j} b_{t-s} r_s $. Using  H\"older's and Rosenthal's \eqref{rosen} inequalities we obtain
\begin{eqnarray}
|\E R_n(t,j)|&\le&\E^{2/3} |r_{t-j}|^3  \E^{1/3}  \big| \sum_{-n < s < t-j} b_{t-s} r_s \big|^3 \nn \\
&\le&\E|r_0|^3 K_3 \big(   \sum_{-n < s < t-j} b^2_{t-s}\big)^{3/2}\  \le \ C.  \label{L2}
\end{eqnarray}
Hence,
\begin{eqnarray}\label{L3}
|h_j|&\le&C|b_j| + \sum_{0<i<j} b_{j-i}^2 |h_i|  \ \le C \big(|b_j| + \sum_{0< i < j} \phi_{j-i} |b_{i}|\big)
\end{eqnarray}
where the first inequality in \eqref{L3}  follows from \eqref{L1} and \eqref{L2} and the second inequality in \eqref{L3}
by iterating the first one
with $\phi_j $ as in \eqref{phi}. Since $\sum_{j=1}^\infty \phi_j = \sum_{k=1}^\infty B^{2k} = 1/(1-B^2)$,
from the second inequality in  \eqref{L3} we obtain $\| h\| \le C B/(1 - B^2)  < \infty. $
The last fact implies $\E R_n(t,j)=  \sum_{i=1}^{n+t-j} h_i b_{i+j} \to \sum_{i=1}^{\infty} h_i b_{i+j} $.
From \eqref{L1} we obtain that the leverage function $h \in L^2(\Z^+) $ is a solution of the linear equation:
\begin{eqnarray}
h_j&=&2a b_j m_2 + \sum_{0<i<j}b_i^2 h_{j-i} +
2b_j \sum_{i>0} b_{i+j} h_{i}, \qquad j =1,2, \cdots.  \label{heq1}
\end{eqnarray}
From Minkowski's inequality, we get  $\sum_{j>0}  (\sum_{0<i<j}b_i^2 h_{j-i})^2 \le B^4 \|h\|^2, $
$\sum_{j>0} \big(b_j \sum_{i>0} b_{i+j} h_{i}\big)^2 $  $ \le  B^4 \|h\|^2 $ and then \eqref{heq1}
implies that $\|h\|\le 2|a|m_2 B + 3 B^2 \|h\| $, or
\begin{eqnarray}
\|h\|&\le&\frac{2|a|m_2 B}{1-3B^2}   \label{heq2}
\end{eqnarray}
provided $B^2 < 1/3$. \\
Let us prove the statements (i) and (ii) of Theorem \ref{lev} for $k=1$. From \eqref{heq1} it follows that
\begin{eqnarray*}
h_1&=&2am_2 b_1+2b_1 \sum_{u=1}^\infty h_u b_{1+u}\ =\  2 b_1(am_2 +\sum_{u=1}^\infty h_u b_{1+u})
\end{eqnarray*}
Since $|\sum_{u=1}^\infty h_u b_{1+u}|\le \|h\| B$, we have
${\rm sgn}(h_1)={\rm sgn}(b_1a)$ provided $\|h\|B<|a|m_2$ holds.
The last relation follows from \eqref{heq2} and $B^2 < 1/5$; indeed,
$$
\|h\|B \le \frac{2|a|m_2 B^2}{1-3B^2}\le |a|m_2.
$$
This proves (i) and (ii) for $k=1$. \\
The general case $k\ge 1$ follows similarly by induction on $k$. Indeed, from \eqref{heq1} we have that
\begin{eqnarray*}
h_k&=&2b_k( am_2+\sum_{u=1}^\infty h_u b_{k+u})+\sum_{j=1}^{k-1}b^2_{k-j}h_j.
\end{eqnarray*}
Assume $h_1,\cdots,h_{k-1}<0 $, then the second term $\sum_{j=1}^{k-1}b^2_{k-j}h_{j}<0$. Moreover,
\begin{equation*}
|\sum_{u=1}^\infty h_u b_{k+u}|\le \|h\|B<|a|m_2
 \end{equation*}
implying that the sign of the first term is the same as ${\rm sgn}(ab_k)$. \\
Theorem \ref{lev} is proved. \hfill $\Box$

\section*{Acknowledgements}

The authors are grateful to two anonymous referees for useful comments.
This research was partly supported by grant no.~MIP-063/2013 from the Research Council of Lithuania.

\end{document}